\newread\epsffilein    
\newif\ifepsffileok    
\newif\ifepsfbbfound   
\newif\ifepsfverbose   
\newdimen\epsfxsize    
\newdimen\epsfysize    
\newdimen\epsftsize    
\newdimen\epsfrsize    
\newdimen\epsftmp      
\newdimen\pspoints     
\pspoints=1bp          
\epsfxsize=0pt         
\epsfysize=0pt         
\def\epsfbox#1{\global\def\epsfllx{72}\global\def\epsflly{72}%
   \global\def\epsfurx{540}\global\def\epsfury{720}%
   \def\lbracket{[}\def\testit{#1}\ifx\testit\lbracket
   \let\next=\epsfgetlitbb\else\let\next=\epsfnormal\fi\next{#1}}%
\def\epsfgetlitbb#1#2 #3 #4 #5]#6{\epsfgrab #2 #3 #4 #5 .\\%
   \epsfsetgraph{#6}}%
\def\epsfnormal#1{\epsfgetbb{#1}\epsfsetgraph{#1}}%
\def\epsfgetbb#1{%
%
%
\openin\epsffilein=#1
\ifeof\epsffilein\errmessage{I couldn't open #1, will ignore it}\else
%
%
   {\epsffileoktrue \chardef\other=12
    \def\do##1{\catcode`##1=\other}\dospecials \catcode`\ =10
    \loop
       \read\epsffilein to \epsffileline
       \ifeof\epsffilein\epsffileokfalse\else
%
%
          \expandafter\epsfaux\epsffileline:. \\%
       \fi
   \ifepsffileok\repeat
   \ifepsfbbfound\else
    \ifepsfverbose\message{No bounding box comment in #1; using defaults}\fi\fi
   }\closein\epsffilein\fi}%
%
%
\def\epsfsetgraph#1{%
   \epsfrsize=\epsfury\pspoints
   \advance\epsfrsize by-\epsflly\pspoints
   \epsftsize=\epsfurx\pspoints
   \advance\epsftsize by-\epsfllx\pspoints
%
%
   \epsfsize\epsftsize\epsfrsize
   \ifnum\epsfxsize=0 \ifnum\epsfysize=0
      \epsfxsize=\epsftsize \epsfysize=\epsfrsize
%
%
     \else\epsftmp=\epsftsize \divide\epsftmp\epsfrsize
       \epsfxsize=\epsfysize \multiply\epsfxsize\epsftmp
       \multiply\epsftmp\epsfrsize \advance\epsftsize-\epsftmp
       \epsftmp=\epsfysize
       \loop \advance\epsftsize\epsftsize \divide\epsftmp 2
       \ifnum\epsftmp>0
          \ifnum\epsftsize<\epsfrsize\else
             \advance\epsftsize-\epsfrsize \advance\epsfxsize\epsftmp \fi
       \repeat
     \fi
   \else\epsftmp=\epsfrsize \divide\epsftmp\epsftsize
     \epsfysize=\epsfxsize \multiply\epsfysize\epsftmp   
     \multiply\epsftmp\epsftsize \advance\epsfrsize-\epsftmp
     \epsftmp=\epsfxsize
     \loop \advance\epsfrsize\epsfrsize \divide\epsftmp 2
     \ifnum\epsftmp>0
        \ifnum\epsfrsize<\epsftsize\else
           \advance\epsfrsize-\epsftsize \advance\epsfysize\epsftmp \fi
     \repeat     
   \fi
%
%
   \ifepsfverbose\message{#1: width=\the\epsfxsize, height=\the\epsfysize}\fi
   \epsftmp=10\epsfxsize \divide\epsftmp\pspoints
   \vbox to\epsfysize{\vfil\hbox to\epsfxsize{%
      \special{illustration #1 scaled \number\epsfscale}
      \hfil}}%
\epsfxsize=0pt\epsfysize=0pt\epsfscale=1000 }%

%
%
{\catcode`\%=12 \global\let\epsfpercent=
%
%
\long\def\epsfaux#1#2:#3\\{\ifx#1\epsfpercent
   \def\testit{#2}\ifx\testit\epsfbblit
      \epsfgrab #3 . . . \\%
      \epsffileokfalse
      \global\epsfbbfoundtrue
   \fi\else\ifx#1\par\else\epsffileokfalse\fi\fi}%
%
%
\def\epsfgrab #1 #2 #3 #4 #5\\{%
   \global\def\epsfllx{#1}\ifx\epsfllx\empty
      \epsfgrab #2 #3 #4 #5 .\\\else
   \global\def\epsflly{#2}%
   \global\def\epsfurx{#3}\global\def\epsfury{#4}\fi}%
%
%
%
%

\newcount\epsfscale    
\newdimen\epsftmpp     
\newdimen\epsftmppp    
\newdimen\epsfM        
\newdimen\sppoints     
\epsfscale=1000        
\sppoints=1000sp       
\epsfM=1000\sppoints
%
\def\computescale#1#2{%
  \epsftmpp=#1 \epsftmppp=#2
  \epsftmp=\epsftmpp \divide\epsftmp\epsftmppp  
  \epsfscale=\epsfM \multiply\epsfscale\epsftmp 
  \multiply\epsftmp\epsftmppp                   
  \advance\epsftmpp-\epsftmp                    
  \epsftmp=\epsfM                               
  \loop \advance\epsftmpp\epsftmpp              
    \divide\epsftmp 2                           
    \ifnum\epsftmp>0
      \ifnum\epsftmpp<\epsftmppp\else           
        \advance\epsftmpp-\epsftmppp            
        \advance\epsfscale\epsftmp \fi          
  \repeat
  \divide\epsfscale\sppoints}
\def\epsfsize#1#2{%
  \ifnum\epsfscale=1000
    \ifnum\epsfxsize=0
      \ifnum\epsfysize=0
      \else \computescale{\epsfysize}{#2}
      \fi
    \else \computescale{\epsfxsize}{#1}
    \fi
  \else
    \epsfxsize=#1
    \divide\epsfxsize by 1000 \multiply\epsfxsize by \epsfscale
  \fi}

\font\cmdunhXII=cmdunh10	at 12pt
\font\sc=cmcsc10								
\font\scXII=cmcsc10 scaled\magstep1

\font\bfVIII=cmbx8

\font\itXIII=cmti10 at 13pt
\font\VIII=cmr8
\font\IX=cmr9
\font\bfIX=cmbx9
\font\go=eufm10

\def\bg{\bigskip\goodbreak}
\def\bgni{\bigskip\goodbreak\noindent}
\def\mg{\medskip\goodbreak}
\def\mgni{\medskip\goodbreak\noindent}

\def\sgni{\smallskip\goodbreak\noindent}
\def\snni{\smallskip\nobreak\noindent}

\def\mnni{\medskip\nobreak\noindent}
\def\bnni{\bigskip\nobreak\noindent}
\let\ni=\noindent
\def\ie{{\it i.e.}}

\newcount\SecNo
\SecNo=0
\newcount\PropNo
\newcount\EqNo

\newif\ifDisplaylines\Displaylinesfalse
{\catcode`\@=11
\gdef\Displaylines#1{\Displaylinestrue\displ@y\halign{
\hbox to\displaywidth{$\@lign\hfil\displaystyle##\hfil$}%
&\llap{$##$}\crcr#1\crcr}}
\gdef\Eqalign#1{\Displaylinestrue\displ@y \tabskip=\centering
\halign to \displaywidth{\hfil$\@lign\displaystyle{##}$\tabskip=0pt
&$\@lign\displaystyle{{}##}$\hfil\tabskip=\centering
&\llap{$\@lign##$}\tabskip=0pt\crcr#1\crcr}}}

\def\PaperNew{\magnification=1100\tolerance=3000
\hsize=14cm\hoffset0.33cm
\vsize=19cm\voffset1.5cm}

\PaperNew

\catcode`\Ç=\active
{\catcode`\@=11
\gdef\Ç#1#2È{\ifcat\noexpand#1\noexpand\%
\global\alloc@4\box\chardef\insc@unt{#1#2}\fi%
\global\setbox#1#2\box0}}
\defÇ#1È#2{\leavevmode\copy#1\ifcat#2a\ \fi#2}
        
\def\EDef{\bigbreak\ni{\bf Definition.}}
\def\EEx{\bigbreak\ni{\bf Example.}}
\def\ERem{\bigbreak\ni{\bf Remark.}}

\def\ECor{\global\advance\PropNo by 1
\setbox0=\hbox{\the\SecNo.\the\PropNo}%
\bigbreak\ni{\bf Corollary
\the\SecNo.\the\PropNo.}}

\def\EEx{\global\advance\PropNo by 1
\setbox0=\hbox{\the\SecNo.\the\PropNo}%
\bigbreak\ni{\bf Example
\the\SecNo.\the\PropNo.}}

\def\ELem{\global\advance\PropNo by 1
\setbox0=\hbox{\the\SecNo.\the\PropNo}%
\bigbreak\ni{\bf Lemma
\the\SecNo.\the\PropNo.}}

\def\Prop{\global\advance\PropNo by 1
\setbox0=\hbox{\the\SecNo.\the\PropNo}%
\bigbreak\ni{\bf Proposition
\the\SecNo.\the\PropNo.}}

\def\Conj{\global\advance\PropNo by 1
\setbox0=\hbox{\the\SecNo.\the\PropNo}%
\bigbreak\ni{\bf Conjecture
\the\SecNo.\the\PropNo.}}

\def\Ques{\global\advance\PropNo by 1
\setbox0=\hbox{\the\SecNo.\the\PropNo}%
\bigbreak\ni{\bf Question
\the\SecNo.\the\PropNo.}}

\def\EThm{\global\advance\PropNo by 1
\setbox0=\hbox{\the\SecNo.\the\PropNo}%
\bigbreak\ni{\bf Theorem
\the\SecNo.\the\PropNo.}}

\def\EPr{\bnni{\it Proof. }}
\def\square{\hbox{$\sqcap\mkern-12mu\sqcup$}}
\def\EndPr\par{{\unskip\nobreak\hfil\penalty50%
\hskip1em\hbox{}\nobreak\square
\parfillskip=0pt\finalhyphendemerits=0\par}}

\def\Eq{\global\advance\EqNo by 1
\global\setbox0=\hbox{\rm(\the\SecNo.\the\EqNo)}
\ifDisplaylines&\copy0\else\eqno{\copy0}\fi}

\newif\iffound
\def\rank#1#2{{\foundfalse\count0=0 \getrank#1#2\end
\iffound\number\count0\else0\fi}} 
\def\getrank#1#2{\ifx#2\end\def\next#1{\relax}%
\else\iffound\else\advance\count0 by1\fi\let\next=\getrank
\ifx#2#1\foundtrue\fi\fi\next#1}

\def\List#1{\def\\##1{\def##1{##1}\rank{##1}{#1}}}

\def\Ref#1; #2; #3; #4\par{\sgni\item{ [#1]}{\sc #2}, {\sl #3}, #4\par}
\def\Reff#1; #2; #3; #4; #5; #6; #7\par{\sgni\item{ [#1]}{\sc #2}, {\sl
#3}, #4 {\bf #5} (#6) #7\par}
						
\newdimen\Dwidth
\newdimen\Dheight
\newdimen\Ddepth

\def\Translate(#1, #2){\dimen105=-#2mm\vglue\dimen105\hskip#1mm}

\def\notdiv{\hbox{$\not |$}}

\def\e{\varepsilon}
\def\scr{\scriptscriptstyle}
\def\dL{\mathord{/}}
\def\dR{\mathord{\setminus}}

\def\gL{\mathbin{\scriptstyle\wedge}}
\def\gR{\mathbin{\widetilde{\scriptstyle\wedge}}}

\def\jL{\mathbin{\widetilde{\scriptstyle\vee}}}
\def\jR{\mathbin{{\scriptstyle\vee}}}

\def\zeq{\buildrel \raise -2pt
				\hbox{$\scr\D$}\over\sim}

\def\dbar{\mid\!\mid}
\def\lg#1{\dbar\!\!#1\!\!\dbar}

\def\D{\Delta}
\let\oo = \infty
\def\resp{{\it resp. }}
\def\s{\sigma}

\def\swedge{{\scriptstyle \wedge}}

\def\sswedge{{\scriptscriptstyle \wedge}}
\def\Permut#1#2{({\raise-0.1pt
				\hbox{$\scriptstyle#1$}
				\atop{\raise1.7pt
				\hbox{$\scriptstyle#2$}}})}
\def\fL{\widetilde\phi}
\def\fR{\phi}

\def\xx{{x}}
\def\yy{{y}}
\def\zz{{z}}
\def\ttt{{t}}

\parindent=12pt

\def\Sec#1\par{\vskip0pt
plus.2\vsize\penalty-250\vskip0pt
plus -.2\vsize
\advance\SecNo by
1\global\PropNo=0\global\EqNo=0
\vskip30pt\centerline{\scXII\the\SecNo.
#1}}

\def\a{\alpha}
\def\b{\beta}

\def\[{\mathrel{\hbox{$[\![$}}}
\def\]{\mathrel{\hbox{$]\!]$}}}

\List{\Adj\Art\Bes\BDM\Bir\BKL\Bou\BMR\BrS\BuZ\ClP
\Dff\Dhh\Dgk\Djj\Dfx\Del\Eps\Gar\HKo\Hig\Hol\LyS
\MKS\Pia\Pib\Pic\Pid\Rem\Rol}

\def\FigCK{1}
\def\FigAutoKnuth{2}

\def\FigArtFourSix{4}
\def\FigRelation{5}
\def\FigWirtThreeTwo{6}
\def\FigWirtThreeFour{7}
\def\FigTypeII{8}
\def\FigWirtFourSix{9}

\def\Chiral{{M_\chi}}
\def\Knuth{{M_\bullet}}

\def\atl#1{A_{\scr#1}}
\def\wir#1{W_{\scr#1}}

\SecNo=0

\bg\centerline{\cmdunhXII AUTOMATIC STRUCTURES FOR
TORUS LINK GROUPS}
\bg\bg\centerline{{\scXII Matthieu
PICANTIN}}

\bg\centerline{ SDAD FRE 2271,
D\'epartement de Math\'ematiques,}

\centerline{Universit\'e de Caen,
Campus II, BP 5186, 14~032 Caen,
France}

\centerline{E-mail: {\tt
picantin@math.unicaen.fr}}

\vskip 0,7 true cm  

\parindent=25pt

{\narrower\ni{\bfVIII Abstract. }{\VIII A general result
of Epstein and Thurston implies that all link groups are
automatic, but the proof provides no explicit automaton.
Here we show that the groups of all torus links are groups
of fractions of so-called Garside monoids, i.e.,
roughly speaking, monoids with a good theory of
divisibility, which allows us to reprove that those
groups are automatic, but, in addition, gives a
completely explicit description of the involved
automata, thus partially answering a question
of~D.~F.~Holt.

\mgni Key words: knot and link groups;
presentation of groups; automatic groups; Artin
groups.

\mgni MSC 2000: 20F05, 57M25, 57M27, 20F10, 20F36,
20M35.\par}}

\parindent=12pt

\def\SecBack{{1}}
\def\SecPres{{2}}
\def\SecGaus{{3}}
\def\SecWirt{{4}}

\vskip30pt\centerline{\scXII Introduction} 

\bnni Epstein and Thurston showed in~[\\\Eps, Chapter XII] that the fundamental group of every
Haken 3-manifold is biautomatic, except if it carries either nilgeometry
or solvgeometry. As a corollary, every link group is biautomatic. However, the proof is not
constructive, in the sense that it provides no explicit automaton.
Here, we propose a combinatorial group theory point of view to
adjoint to this geometric group theory result. The aim of
this paper is to show that all the torus link groups are the
groups of fractions of monoids that admit a nice theory of
divisibility, like spherical Artin or Birman-Ko-Lee monoids~:
all the torus link groups belong so to the wide class of
so-called {\it Garside groups}. Dehornoy stated
in~[\\\Dgk] that every Garside group is biautomatic, by
proving that each Garside monoid of which it is the
group of fractions provides a simple explicit
biautomatic structure. Therefore, our result enables us to
construct an explicit finite state automaton computing normal
forms for each torus link group. 

\mgni This paper is organized as follows. In Section~\SecBack, we gather
earlier results of~[\\\Dgk] and~[\\\Dfx] about~Garside
monoids and groups.  In Section~\SecPres, we compute a first well-fitted
(monoid) presentation for every torus link group, then we complete the
latter to obtain a so-called {\it complemented} presentation.
In Section~\SecGaus, using the previous results, we show that
every torus link group is a Garside group, and,
therefore, that it admits a biautomatic structure.
Several examples are studied. Finally, in Section~{\SecWirt},
we discuss the existence of other automatic structures for
these groups.

\mnni The author wishes to thank~D.~B.~A.~Epstein,
S.~H.~Gersten and~D.~F.~Holt for valuable discussions.

\Sec Background from Garside groups

\bnni In this section, we list some
basic properties of Garside monoids and groups, and summarize
Dehornoy's results about biautomaticity of Garside
groups. Finally, we recall an effective criterion for
recognizing a Garside monoid from a given monoid
presentation. For all the results quoted here, we refer the
reader to~[\\\Dfx], [\\\Dgk] and~[\\\Dhh].

\mg\bgni{\itXIII Main definitions and basic properties}

\mnni Assume that~$M$ is a
monoid. We say that~$M$ is {\it conical}
if~$1$ is the only invertible
element in~$M$. For~$a,b$
in~$M$, we say that~$b$ is a left
divisor of~$a$---or that~$a$ is a
right multiple
of~$b$---if~$a=bd$ holds for
some~$d$ in~$M$. An
element~$c$ is a right lower
common multiple---or a right
lcm---of~$a$ and~$b$ if it is a
right multiple of both~$a$
and~$b$, and every common
right multiple of~$a$ and~$b$ is
a right multiple of~$c$. Right
divisor, left multiple, and left
lcm are defined symmetrically.
For~$a,b$ in~$M$, we say
that~$b$ {\it divides}~$a$---or
that~$b$ is a divisor
of~$a$---if~$a=cbd$ holds for
some~$c,d$ in~$M$. 

If~$c$, $c'$ are two right lcm's
of~$a$ and~$b$, necessarily~$c$
is a left divisor of~$c'$, and~$c'$
is a left divisor of~$c$. If we
assume~$M$ to be conical and
cancellative, we have~$c=c'$. In
this case, the unique right lcm
of~$a$ and~$b$ is denoted by~$a
\jR b$. If~$a\jR b$ exists,
and~$M$ is left cancellative,
there exists a unique
element~$c$ satisfying~$a
\jR b=ac$. This element is
denoted by~$a\dR b$. We define
the {\it left lcm}~$\jL$ and the
left operation~$\dL$
symmetrically.  In particular,
we have$$a\jR b=a(a\dR
b)=b(b\dR a),
\hbox{\quad  and\quad} a\jL
b=(b\dL a)a=(a\dL b)b.$$Cancellativity
and conicity implies that
left and right divisibility are
order relations.

\EDef\Ç\GaussianMonoidÈ[\\\Dgk]
A monoid~$M$ is said to be {\it
Garside} if $M$ is conical and
cancellative, every pair of
elements in~$M$ admits a left
lcm and a right lcm, and~$M$ admits a {\it Garside element},
defined to be an element whose left and right divisors
coincide and generate~$M$.

\bgni By~[\\\BrS], all spherical
Artin monoids are Garside monoids. The braid
monoids of the complex
reflection
groups~$G_7,G_{11},G_{12},G_{13},
G_{15},G_{19}$ and~$G_{22}$
given in~[\\\BMR] (see~[\\\Dfx][\\\Pic]), \hbox{Garside's
hypercube monoids ([\\\Gar][\\\Pic]),} the
Birman-Ko-Lee monoids for spherical Artin
groups ([\\\BKL][\\\BDM][\\\Pic][\\\Bes])
are also Garside monoids.

\EEx\Ç\ExChiralÈ{Let us consider the
monoid~$\Chiral$ with presentation
$$\langle~x,y,z: xzxy=yzx^2~,~ yzx^2z=zxyzx~,~
zxyzx=xzxyz~\rangle.$$ The
monoid~$\Chiral$ is a typical example of a
Garside monoid, and, in addition,
$\Chiral$ has the distinguishing feature to
be not antiautomorphic, contrary to all
spherical Artin monoids.}

\bnni If~$M$ is a Garside monoid, then~$M$
satisfies Ore's
conditions~[\\\ClP], and it
embeds in a group of right
fractions, and, symmetrically, in
a group of left fractions. In this
case, by construction, every
right fraction~$ab^{-1}$
with~$a,b$ in~$M$ can be
expressed as a left
fraction~$c^{-1}d$, and
conversely. Therefore, the two
groups coincide, and  there is no
ambiguity in speaking of {\it
the} group of fractions of a Garside monoid.

\EDef\Ç\GaussianGroupÈ{A
group~$G$ is a {\it Garside group} if there exists a
Garside monoid of
which~$G$ is the group of
fractions.}

\ELem\Ç\CalculousÈ[\\\Dgk]
{\sl Assume that~$M$ is a
Garside monoid. Then the
following identities holds
in~$M$:$$\Displaylines{(ab)\jR(ac)=a(b\jR
c),\Eq\Ç\IdentityIÈ\cr
c\dR(ab)=(c\dR a)((a\dR c)\dR
b),\hbox{\qquad} (ab)\dR c
=b\dR(a\dR
c),\Eq\Ç\IdentityIIÈ\cr (a\jR
b)\dR c =(a\dR b)\dR(a\dR c)
=(b\dR a)\dR(b\dR
c),\hbox{\qquad} c\dR(a\jR
b)=(c\dR a)\jR(c\dR
b).\Eq\Ç\IdentityIIIÈ}$$

}\ELem\Ç\thinAtomicÈ[\\\Dfx]
{\sl Assume that~$M$ is a Garside monoid. Then the
following equivalent assertions
hold~:

\snni (i) There exists a
mapping~$\mu$ from~$M$ into
the integers
satisfying~$\mu(a)>0$ for
every~${a\not=1}$ in~$M$,
and
satisfying~$\mu(ab)\geq\mu(a)+\mu(b)$
for every~$a,b$ in~$M$; 

\snni (ii) For every set~$X$ that
generates~$M$ and for
every~$a$ in~$M$, the lengths of
the decompositions of~$a$ as
products of elements in~$X$
have a finite upper bound.}

\EDef\Ç\NormÈ[\\\Dfx] {A monoid is said to be {\it
atomic} if it satisfies
the equivalent conditions of
Lemma~Ç\thinAtomicÈ . An {\it atom} is defined to
be a non trivial element~$a$
such that~$a=bc$
implies~$b=1$ or~$c=1$.
The {\it norm} function~\hbox{$\lg{.}$} of
an atomic monoid~$M$ is defined such that, for
every~$a$ in~$M$,
$\lg{a}$ is the upper bound of the lengths of the
decompositions of~$a$ as products of atoms.}

\EEx\Ç\ExKnuthÈ{The monoid~$\Knuth$ defined by
the presentation$$\langle~x,y: xyxyx=yy~\rangle$$
is an other example of a Garside monoid, which,
as for it, admits no additive norm, \ie, no norm~$\nu$
satisfying~$\nu(ab)=\nu(a)+\nu(b)$.

\bgni By the previous lemma,
every element in a Garside monoid has
finitely many left divisors, only then,
for every pair of elements~$(a,
b)$, the common left divisors
of~$a$ and~$b$ admit a right
lcm, which is therefore the {\it 
left gcd} of~$a$ and~$b$. This
left gcd will be denoted
by~$a\gL b$. We define  the {\it
right gcd}~$\gR$ symmetrically. 

\ELem\Ç\ClosureÈ[\\\Dgk] {\sl
Assume that~$M$ is a Garside monoid. Then  it
admits a finite generating
subset that is closed
under~$\dR,\dL,\jR,\gL,\jL$
and~$\gR$.}

\bgni Every Garside monoid
admits a finite set of atoms, and
this set is the minimal
generating set~[\\\Dfx]. The
hypothesis that there exists a
finite generating subset that is
closed under~$\dR$ implies that
the closure of the atoms
under~$\dR$ is finite---its
elements are called {\it right
primitive elements}. In
particular, the closure of the
atoms under~$\dR$ and~$\jR$ is
finite---its elements are called
{\it simple elements}, and their
right lcm is denoted by~$\D$. It
turns out that the set of the
simple elements is also the
closure of atoms under~$\dL$
and~$\jL$. So, the element~$\D$
is both the right and the left lcm
of the simple elements, and it is
called {\it the Garside element}
of the monoid. If~$M$ is a Garside monoid and~$S$ is the
set of simple elements in~$M$,
then~$(S,\gL,\jR,1,\D)$ is a
finite lattice, which completely determines~$M$.

\EEx\Ç\LatticeÈ{The lattices of
simple elements in~$\Chiral$
(Example~Ç\ExChiralÈ) and in~$\Knuth$
(Example~Ç\ExKnuthÈ) are displayed
in~Figure~\FigCK. }

\def\newlabel <#1,#2> #3 {\smash{\rlap{\kern #1mm \raise #2mm\hbox{#3}}}}

\midinsert
\centerline{
\newlabel <23, -3.3> $1$
\newlabel <-2, 12> $\xx$
\newlabel <19.5, 12> $\yy$
\newlabel <48, 12> $\zz$
\newlabel <22.8, 73> $\D$
\newlabel <89, -3.3> $1$
\newlabel <73.5, 12> $\xx$
\newlabel <104.5, 12> $\yy$
\newlabel <88.5, 73> $\D$
\epsfscale=450
\epsfbox{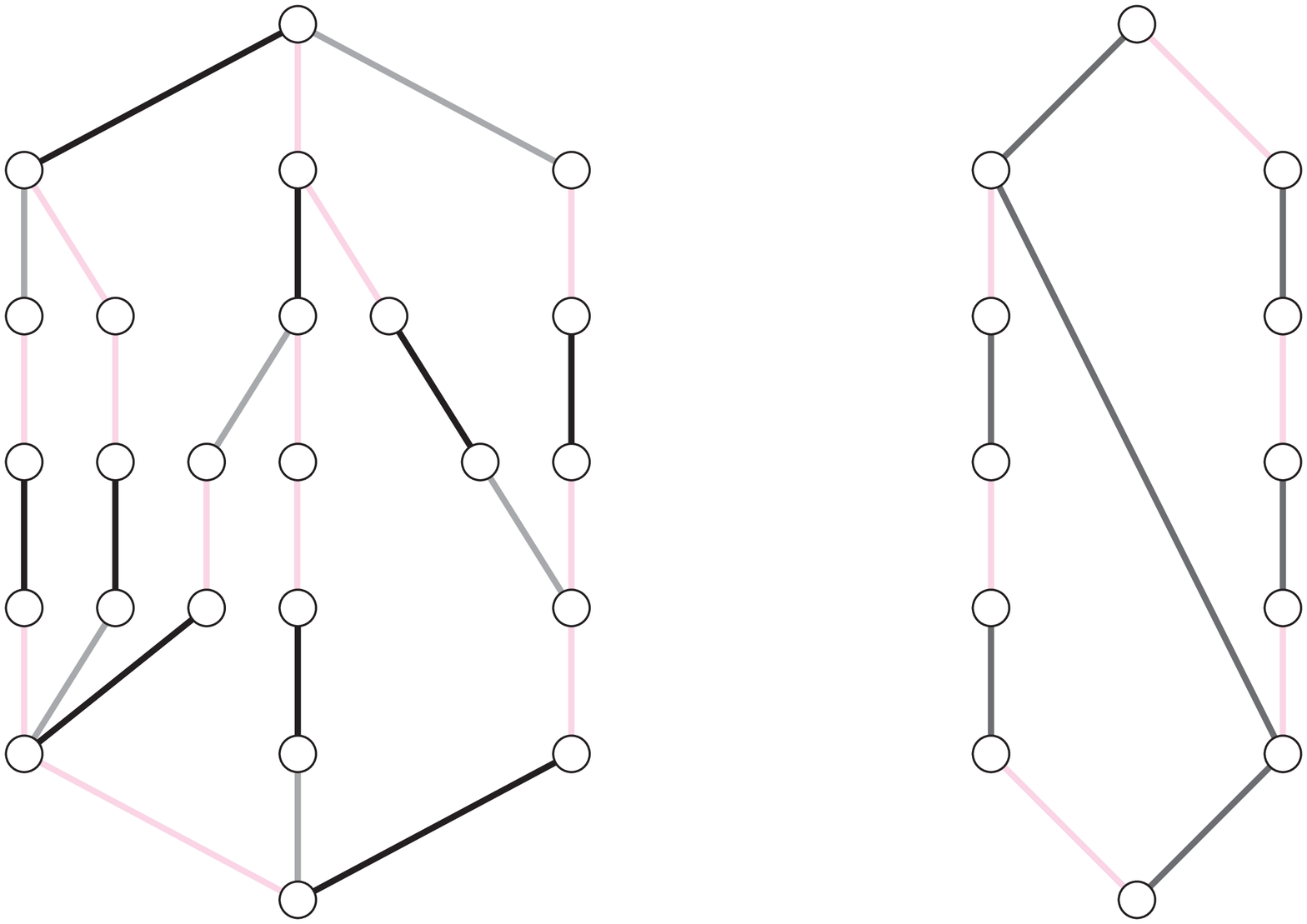}
}
\vglue 3mm
\centerline{\IX{\bfIX Figure~\FigCK.} The
lattices of simples in~$\Chiral$ (on the left) and
in~$\Knuth$ (on the right).}
\endinsert

\Prop\Ç\DehÈ [\\\Dgk] {\sl
Assume that~$M$ is a Garside monoid,
$S$ is the set of its simple
elements, and~$\D$ is its
Garside element.

\snni(i) Let~$k$ be a
nonnegative integer. Then,
$S^k$ is both the set of all left
divisors of~$\D^k$ and the set of
all right divisors of~$\D^k$.

\snni(ii) The
functions~$a\mapsto
(a\dR\D)\dR\D$
and~$a\mapsto\D\dL (\D\dL a)$
from~$S$ into itself extend into
automorphisms~$\fR$ and~$\fL$
of~$M$ that map~$S^k$ into
itself for every~$k$, and the
equalities$$a\D =
\D\fR(a),
\hbox{\quad and
\quad} \D  a = \fL(a)\D$$hold for
every~$a$ in~$M$.}

\mg\bgni {\itXIII Normal forms and automatic structures}

\mnni If~$S$ generates the monoid~$M$, every element of~$M$ can be
decomposed as a product of elements of~$S$. The classical idea for
obtaining a unique distinguished decomposition is to push the elements to
the right (or to the left), \ie, to require the last element of the
decomposition to be the maximal possible one. So a
$\dR$-closed family is what is needed for this process to
lead to a unique decomposition
(see~[\\\Dfx][\\\Dgk][\\\Dhh][\\\Pia]).

\EDef\Ç\GreedyÈAssume that~$M$ is a Garside monoid, and~$G$ is its
group of fractions. For~$c$ in~$G$, the {\it (right fractional) normal
form} of~$c$ is defined to be the unique decomposition$$a_p\cdots
a_1b_1^{-1}\cdots b_q^{-1},$$where~$a_p,\ldots,b_q$ are simple elements
in~$M$ all distinct from~1, and we have~$a_1\gR
b_1=1$, $a_i\gR(\D\dL a_{i-1})=1$ for~$2\leq i\leq p$,
and~$b_i\gR(\D\dL b_{i-1})=1$ for~$2\leq i\leq q$.

\bgni The previous normal forms are associated with a biautomatic
structure. More than an existence result, simple
automata that compute normal forms can be constructed~:

\Prop\Ç\AutomataÈ[\\\Dgk] {\sl Assume that $M$ is a Garside
monoid. Let~$A$ be the set of atoms, $S$ be the set of simples and~$\D$
be the Garside element. Let define~$T:S\times A\rightarrow S$
by~$T(s,\xx)=((\D\dL(\D\dL\xx))\dR(\D\dL s))\dR\D$. Then, for
every word~$u$ on~$A$, we have~$\bar{u}\gR\D=T(1,u)$, \ie, the
result of reading~$u$ by the automaton~$(T,1)$ is the last
term of the right normal~form~of~$\bar{u}$.}

\ECor\Ç\AutomaticÈ[\\\Dgk] {\sl Every Garside group admits a
biautomatic structure.}

\EEx{ From the lattice of simples in~$\Knuth$ (Figure~\FigCK),
one can easily construct a transducer---\ie, a finite state
automaton with output---computing the right normal form of
every element in~$\Knuth$~: see Figure~{\FigAutoKnuth}. The
initial state is represented with a big black arrow~: the
state~1. Clear (\resp dark) arrows represent the reading of
the letter~$\xx$ (\resp the letter~$\yy$). During the reading
of a word~$w$, one concatenates the labels (eventually
empty) of crossed arrows. At the end of the reading of~$w$,
the ambient state~$s$ is the first simple of the right
normal form~$N(w)$ of~$w$ and the word~$w'$ obtained by
concatenating the various labels is the word that remains to
normalize~: we have~$N(w)=N(w')\cdot s$.}

\midinsert
\centerline{
\epsfscale=500
\epsfbox{PieFigure2.eps}
}
\vglue 3mm
\centerline{\IX{\bfIX Figure~\FigCK.} A
transducer computing the right normal form in~$\Knuth$.}
\endinsert

\mg\bgni {\itXIII
Recognizing Garside monoids}

\mnni We conclude this section by recalling how to effectively recognize
Garside monoids from a presentation. In fact, a correspondence can
be established between Garside monoids and presentations of a
particular syntactic form. The connection relies on a combinatorial technique
called {\it word redressing} (word reversing in [\\\Dhh]), and it leads to two fundamental
applications. On the one hand, it allows to effectively decide whether a
given presentation defines a Garside monoid. On the
other hand, given a Garside monoid presentation, it
allows one to compute the main operations of the monoid.

\bgni For a nonempty set~$A$---usually called an {\it
alphabet} in this context---we denote by~$A^*$ the free
monoid with basis~$A$. We use~$\e$ for the empty word.

\EDef\Ç\ComplementÈ[\\\Dgk] {We say that a monoid
presentation~$\langle~A:R~\rangle$ is {\it complemented} if,
for every pair~$(x,y)$ in~$A^2$, there exists exactly one
relation~$xu=yv$ in~$R$~; in this case, $u$ and~$v$ are
denoted by~$x\dR y$ and~$y\dR x$ respectively. The
operation~$\dR$ on~$A$ can be extended into a unique partial
operation~$\dR$ on~$A^*$ satisfying the rules of word
redressing~:
$$\eqalign{(uv)\dR w&=v\dR(u\dR w),\cr
w\dR(uv)&=(w\dR u)((u\dR w)\dR v).}$$
A complemented monoid presentation~$\langle~A:R~\rangle$ is
called {\it coherent} on~$B\in A^*$ if, for every~$(u,v,w)$
in~$B^3$, we have
$$((u\dR v)\dR(u\dR
w))\dR((v\dR u)\dR(v\dR w))=\e,$$and called
{\it normed} if there exists a mapping~$\nu$ of~$A^*$ to the
integers such that, for all~$\xx,\yy$
in~$A$, and all~$u,v$ in~$A^*$, we have~$\nu(\xx v)>\nu(v)$,
$\nu(v\xx)>\nu(v)$
and~$\nu(u\xx(\xx\dR\yy)v)=\nu(u\yy(\yy\dR\xx)v)$.}

\bgni The Garsidity criterion we shall use in the sequel
is~:

\Prop\Ç\CriterionÈ[\\\Dhh] {\sl Assume that~$M$ is a monoid
with a finite coherent normed and complemented presentation.
Assume that~$M$ is cancellative and admits a Garside
element. Then $M$ is a Garside monoid.}

\bgni For alternative Garsidity criteria and details, we
refer to~[\\\Dgk], see also~[\\\Dhh][\\\Djj][\\\Dfx].

\ERem{ We know no such criterion allowing to prove
that a given finite complemented presentation define a (no
Garside) Gaussian monoid,
\ie, a Gaussian monoid without a Garside element. The
coherence condition is certainly necessary, but not
sufficient. Let us mention moreover that we know actually no
example of finitely generated Gaussian monoid or group that is
not Garside~: the braid group~$B_\oo$ is Gaussian and not
Garside, but it is not finitely generated.}

\vfill\eject
\Sec Monoid presentations for torus link groups

\def\sbul{\mathrel{\scr\bullet}}

\bnni In this section, we compute a first well-fitted (monoid)
presentation for every torus link group, and then complete
the latter to obtain a complemented presentation.

\mgni Every oriented link can be realized as a closed braid
$\widehat b$ for some braid~$b$. Following
Artin [\\\Art] (see also~[\\\Bir]), if~$b$ is an $n$-strand braid, then
the link group of~$\widehat b$---\ie, the fundamental
group~$\pi_1(S^3\setminus\widehat b)$ of the complement
of~$\widehat b$---has a presentation of the form
$$\langle~\xx_1,\ldots,\xx_n:\xx_i=\xx_i\sbul b~~\hbox{for}~1\leq i\leq
n~\rangle,$$ where $u\sbul b$ denotes the image of the
element~$u$ under some automorphism~$b$ of the free group
based on~$\{\xx_1,\ldots,\xx_n\}$.

\EDef\Ç\TorusLinkÈ A {\it torus link} is a link that can be
drawn on the surface of a standard torus without
self-intersection~: given integers~$p$ and~$q$, the
$(p,q)$-torus link wraps meridionally around the torus
$p$~times and wraps longitudinally around the torus
$q$~times.

\bgni Applying Artin's theorem to the particular case of torus links, we
obtain that, for~$2\leq p\leq q$, the group of the ($p,q$)-torus link is
$$\langle~\xx_1,\ldots,\xx_p:
\xx_i=\xx_i\sbul(\s_1\cdots\s_{p-1})^q
\hbox{ for }1\leq i\leq p~\rangle,$$ since the $(p,q)$-torus link is the
closure of the braid~$(\s_1\cdots\s_{p-1})^q$, where the $\s_i$'s denote the
Artin generators of~$B_n$.

\ELem\Ç\TorusLinkGroupÈ{\sl For~$2\leq p\leq q$, the ($p,q$)-torus link
group is
$$\langle~\xx_1,\ldots,\xx_p:
\xx_i(\xx_1\cdots\xx_p)^{[{{q-i}\over{p}}]+1}=(\xx_1\cdots\xx_p)^{[{{q-i}\over{p}}]+1}
\xx_{p\lambda_i-q+i}
\hbox{ for }1\leq i\leq p~\rangle.\Eq\Ç\originalÈ$$

}\EPr We consider the action of the braid group~$B_p$ onto the free group
on~$\{\xx_1,\ldots,\xx_p\}$ defined by
$$\xx_i\sbul\s_j=\cases{\xx_i\xx_{i+1}\xx_i^{-1}&for~$i=j$\cr
\xx_{i-1}&for~$i=j+1$\cr\xx_i&otherwise.}$$ Let~$[z]$ denote the integer
part of~$z$. We show by induction on~$q\geq2$ that, for every~$p$
with~$2\leq p\leq q$, we have
$$\xx_i\sbul(\s_1\cdots\s_{p-1})^q
=(\xx_1\cdots\xx_p)^{[{{q-i}\over{p}}]+1}
\xx_{p[{{q-i}\over{p}}]+p-q+i}
(\xx_1\cdots\xx_p)^{-[{{q-i}\over{p}}]-1}$$
for every~$1\leq i\leq p$. For~$p=q=2$, we
find$$\eqalign{\xx_1\sbul\s_1^2&=\xx_1\xx_2\xx_1^{-1}\sbul\s_1
=(\xx_1\xx_2\xx_1^{-1})(\xx_1)(\xx_1\xx_2^{-1}\xx_1^{-1})
=(\xx_1\xx_2)\xx_1(\xx_1\xx_2)^{-1},\cr
\xx_2\sbul\s_1^2&=\xx_1\sbul\s_1
=\xx_1\xx_2\xx_1^{-1}=(\xx_1\xx_2)\xx_2(\xx_1\xx_2)^{-1}.}$$Assume~$q>2$.
By induction hypothesis, we have
$$\eqalign{\xx_i\sbul(&\s_1\cdots\s_{p-1})^q\cr
&=(\xx_1\cdots\xx_p)^{[{{q-i-1}\over{p}}]+1}
\xx_{p[{{q-i-1}\over{p}}]+p-q+i+1}
(\xx_1\cdots\xx_p)^{-[{{q-i-1}\over{p}}]-1}
\sbul(\s_1\cdots\s_{p-1}),}$$
for~$1\leq i<p$. We easily
compute~$(\xx_1\cdots\xx_p)^k\sbul\s_i=(\xx_1\cdots\xx_p)^k$
for~$1\leq i<p$,
hence$$(\xx_1\cdots\xx_p)^k\sbul(\s_1\cdots\s_{p-1})
=(\xx_1\cdots\xx_p)^k$$ for~$k\geq 0$ (also a
direct consequence of the fact that the considered  action is an
automorphism), and$$\xx_i\sbul(\s_1\cdots\s_{p-1}) =\cases{(\xx_1\cdots
\xx_p)\xx_p(\xx_1\cdots\xx_p)^{-1}&for~$i=1$\cr
\xx_{i-1}&for~$1<i\leq p$.}$$
Now, $p[{{q-i-1}\over{p}}]+p-q+i+1=1$
is equivalent
to~$[{{q-i-1}\over{p}}]+1
={{q-i}\over{p}}=[{{q-i}\over{p}}]$.
We deduce
$$\xx_1\sbul (\s_1\cdots\s_{p-1})^q
=(\xx_1\cdots\xx_p)^{[{{q-i}\over{p}}]}
(\xx_1\cdots
\xx_p)\xx_p(\xx_1\cdots\xx_p)^{-1}
(\xx_1\cdots\xx_p)^{-[{{q-i}\over{p}}]},$$
and$$\xx_i\sbul(\s_1\cdots\s_{p-1})^q
=(\xx_1\cdots\xx_p)^{[{{q-i}\over{p}}]+1}
\xx_{p[{{q-i}\over{p}}]+p-q+i}
(\xx_1\cdots\xx_p)^{-[{{q-i}\over{p}}]-1},$$
for~$1<i<p$. Finally, for~$i=p$, a direct computation
gives the required
result$$\xx_p\sbul(\s_1\cdots\s_{p-1})^q
=(\xx_1\cdots\xx_p)^{[{q\over
p}]}\xx_{p[{q\over p}]+p-q}
(\xx_1\cdots\xx_p)^{-[{q\over p}]},$$ which completes
the induction.\EndPr

\bgni Now we have to distinguish two cases in order to
complete the presentation obtained above into a new one, that
will be eligible for applying the Garsidity criterion
of Proposition~Ç\CriterionÈ.

\mgni By Lemma~Ç\TorusLinkGroupÈ, for~$p\geq 2$
and~$k\geq 1$, the ($p,pk$)-torus link group is
$$\langle~\xx_1,\ldots,\xx_p:
\xx_i(\xx_1\cdots\xx_p)^k
=(\xx_1\cdots\xx_p)^k
\xx_i
\hbox{ for }1\leq i\leq p~\rangle,$$
\ie,$$\langle~\xx_1,\ldots,\xx_p:
(\xx_1\cdots\xx_p)^k
=(\xx_2\cdots\xx_p\xx_1)^k
=\ldots=(\xx_p\xx_1\cdots\xx_{p-1})^k~\rangle.
\Eq\Ç\ArtinTorÈ$$
\bgni Let us fix~$2\leq
p\leq q$ and assume that~$p$ does not divide~$q$. Let~$\a,\b$
be the integer part and~the remainder of~$q/p$. Then we
have~$q=p\a+\b$ with~$1\leq
\b\leq p-1$. The ($p,q$)-torus link group is then
$$\eqalign{\langle~\xx_1,\ldots,\xx_p:
(\xx_1\cdots\xx_p)^{\a+1}
&=(\xx_2\cdots\xx_p)(\xx_1\cdots\xx_p)^\a
\xx_{p-\b+1}\cr
&=(\xx_3\cdots\xx_p)(\xx_1\cdots\xx_p)^\a
\xx_{p-\b+1}\xx_{p-\b+2}\cr
&=\ldots\cr
&=(\xx_\b\cdots\xx_p)(\xx_1\cdots\xx_p)^\a
\xx_{p-\b+1}\cdots\xx_{p-1}\cr
(\xx_1\cdots\xx_p)^\a\xx_1
&=\xx_{\b+1}(\xx_1\cdots\xx_p)^\a\cr
(\xx_1\cdots\xx_p)^\a\xx_2
&=\xx_{\b+2}(\xx_1\cdots\xx_p)^\a\cr
&\;\;\vdots\cr
(\xx_1\cdots\xx_p)^\a\xx_{p-\b}
&=\xx_p(\xx_1\cdots\xx_p)^{\a}~\rangle.}
\Eq\Ç\ArtinTorIIÈ$$

\EDef{ Let~$2\leq
p\leq q$. We call {\it $(p,q)$-torus link Artin monoid} the
monoid~$\atl{p,q}$ with
presentation~Ç\ArtinTorÈ for~$p|q$ and with
presentation~Ç\ArtinTorIIÈ for~$p\notdiv q$.}

\bgni By construction, for~$2\leq p\leq q$, the group admitting the
presentation of~$\atl{p,q}$ is the~$(p,q)$-torus link group. Now, we show
that $\atl{p,q}$ admits a complemented presentation. This new
presentation (not unique {\it a priori}) will allow us in the
following section to show that~$\atl{p,q}$ is a Garside monoid,
which, in particular, will tell us that this complemented presentation is
the right one, \ie, the minimal one.

\Prop\Ç\ComplementedÈ{\sl Let~$2\leq
p\leq q$ with~$p\notdiv q$. The monoid~$\atl{p,q}$
admits the complemented
presentation$$\langle~\xx_1,
\ldots,\xx_p~:~\xx_i(\xx_i\dR\xx_j)
=\xx_j(\xx_j\dR\xx_i)\hbox{
for }1\leq i,j\leq p~\rangle,$$
where~$\dR$ is defined by

$$\xx_{s\b+1}\dR\xx_{t\b+1}
=\cases{
w_{\scr(t\a)0(\a-1)}
&for~$t<s$,\cr
w_{\scr(s\a)1(\a-1)}\xx_{(t-s-1)\b+1}
&for~$t>s$,\cr}$$
by$$\xx_{s\b+e}\dR\xx_{t\b+1}
=\cases{
w_{\scr(t\a)0(\a-1)}
&for~$t<s$,\cr
w_{\scr(s\a)e\a}
&for~$s\leq t\leq s+1$,\cr
w_{\scr(s\a)e\a}w_{\scr00(\a-2)}\xx_{(t-s-2)\b+1}&for~$s+1<t$,\cr}$$
and
$$\xx_{t\b+1}\dR\xx_{s\b+e}
=\cases{
w_{\scr(t\a)1(\a-1)}\xx_{(s-t-1)\b+e}
&for~$t<s$,\cr
w_{\scr(t\a)1(s\a-t\a+\a)}
&for~$s\leq t\leq s+1$,\cr
w_{\scr(s\a)0(2\a-1)}
&for~$s+1<t$,\cr}$$
for~$e\not=1$, and by
$$\xx_{s\b+e}\dR\xx_{t\b+f}
=\cases{
w_{\scr(t\a)0(2\a-1)}
&for~$t<s-1$,\cr
w_{\scr(s\a)e\a}
&for~$t=s-1$, or for~$s=t$ with~$e\not=f$,\cr
w_{\scr(s\a)e\a}w_{\scr00(\a-1)}
&for~$t=s+1$,\cr
w_{\scr(s\a)e\a}w_{\scr00(\a-2)}\xx_{(t-s-2)\b+f}
&for~$t>s+1$,\cr}$$
for~$e,f\not=1$, with$$w_{\scr
ze\a}=(\xx_1\cdots\xx_p)^z\xx_{e+1}\cdots\xx_p
(\xx_1\cdots\xx_p)^\a\xx_{1+p-\b}\cdots\xx_{e-1+p-\b},$$
implying~$w_{\scr z0\a}=(\xx_1\cdots\xx_p)^{z+1+\a}$
and~$w_{\scr z1\a}=(\xx_1\cdots\xx_p)^z\xx_2\cdots\xx_p
(\xx_1\cdots\xx_p)^\a$.}

\EPr First, the relations in~Ç\ArtinTorIIÈ are
immediate consequences of those of the presentation above.
Conversely, we have to find a relation of the
form~$\xx_i\cdots=\xx_j\cdots~$ for~$1\leq i,j\leq
p$ using only the relations of the presentation~Ç\ArtinTorIIÈ
for~$\atl{p,q}$, that is, using only the relations
$$\eqalignno{\xx_e\cdots\xx_p(\xx_1\cdots
\xx_p)^\a\xx_{p-\b+1}\cdots\xx_{e-1+p-\b}
&=\xx_f\cdots\xx_p(\xx_1\cdots\xx_p)^\a
\xx_{p-\b+1}\cdots\xx_{f-1+p-\b},&(1)\cr
\xx_{s\b+e}(\xx_1\cdots\xx_p)^{s\a}
&=(\xx_1\cdots\xx_p)^{s\a}\xx_e,&(2)}$$
for~$1\leq e,f\leq \b$. 

\mgni{\it Case 1.} Assume~$i=s\b+e$
and~$j=t\b+1$ with~$t<s$. Then we have
$$\eqalign{\xx_{s\b+e}(\xx_1\cdots\xx_p)^{(t+1)\a}
&=^{{}^{\!\!\!\!\!\!\!\rm\scr(2)}}\!
(\xx_1\cdots\xx_p)^{(t+1)\a}\xx_{(s-t-1)\b+e}\cr
&=(\xx_1\cdots\xx_p)^{t\a}\xx_1\xx_2\cdots\xx_p(\xx_1\cdots\xx_p)^{\a-1}\xx_{(s-t-1)\b+e}\cr
&=^{{}^{\!\!\!\!\!\!\!\rm\scr(1)}}\!
\xx_{t\b+1}(\xx_1\cdots\xx_p)^{t\a}\xx_2\cdots\xx_p(\xx_1\cdots\xx_p)^{\a-1}\xx_{(s-t-1)\b+e}.}$$
This gives both the expression
of~$\xx_{s\b+1}\dR\xx_{t\b+1}$ for~$0\leq
s\not=t$ and the expressions
of~$\xx_{s\b+e}\dR\xx_{t\b+1}$
and~$\xx_{t\b+1}\dR\xx_{s\b+e}$ for~$0\leq s<t$
and~$e\not=1$.

\mgni{\it Case 2.} Assume~$i=s\b+e$
and~$j=t\b+1$ with~$s\leq t\leq s+1$ and~$e\not=1$. Then we have
$$\eqalign{\xx_{s\b+e}&(\xx_1\cdots\xx_p)^{s\a}
\xx_{e+1}\cdots\xx_p(\xx_1\cdots\xx_p)^\a
\xx_{1+p-\b}\cdots\xx_{e-1+p-\b}\cr
&=^{{}^{\!\!\!\!\!\!\!\rm\scr(2)}}\!
(\xx_1\cdots\xx_p)^{s\a}
\xx_e\xx_{e+1}\cdots\xx_p(\xx_1\cdots\xx_p)^\a
\xx_{1+p-\b}\cdots\xx_{e-1+p-\b}\cr
&=^{{}^{\!\!\!\!\!\!\!\rm\scr(1)}}\!
(\xx_1\cdots\xx_p)^{(s+1)\a+1}=(\xx_1\cdots\xx_p)^{t\a}\xx_1\xx_2\cdots\xx_p(\xx_1\cdots\xx_p)^{(s-t+1)\a}\cr
&=^{{}^{\!\!\!\!\!\!\!\rm\scr(2)}}\!
\xx_{t\b+1}(\xx_1\cdots\xx_p)^{t\a}\xx_2\cdots\xx_p(\xx_1\cdots\xx_p)^{(s-t+1)\a}.}$$
This provides the expressions
of~$\xx_{s\b+e}\dR\xx_{t\b+1}$
and~$\xx_{t\b+1}\dR\xx_{s\b+e}$ for~$0\leq s\leq
t\leq s+1$ and~$e\not=1$.

\mgni{\it Case 3.} Assume~$i=s\b+e$
and~$j=t\b+1$ with~$s+1<t$ and~$e\not=1$. Then we have
$$\eqalign{\xx_{s\b+e}&(\xx_1\cdots\xx_p)^{s\a}
\xx_{e+1}\cdots\xx_p(\xx_1\cdots\xx_p)^\a
\xx_{1+p-\b}\cdots\xx_{e-1+p-\b}(\xx_1\cdots\xx_p)^{\a-1}\xx_{(t-s-2)\b+1}\cr
&=^{{}^{\!\!\!\!\!\!\!\rm\scr(2)}}\!
(\xx_1\cdots\xx_p)^{s\a}
\xx_e\xx_{e+1}\cdots\xx_p(\xx_1\cdots\xx_p)^\a
\xx_{1+p-\b}\cdots\xx_{e-1+p-\b}(\xx_1\cdots\xx_p)^{\a-1}\xx_{(t-s-2)\b+1}\cr
&=^{{}^{\!\!\!\!\!\!\!\rm\scr(1)}}\!
(\xx_1\cdots\xx_p)^{(s+2)\a}\xx_{(t-s-2)\b+1}
=^{{}^{\!\!\!\!\!\!\!\rm\scr(2)}}\!
\xx_{t\b+1}(\xx_1\cdots\xx_p)^{(s+2)\a}.}$$
This gives the expressions
of~$\xx_{s\b+e}\dR\xx_{t\b+1}$
and~$\xx_{t\b+1}\dR\xx_{s\b+e}$ for~$0<s+1<t$
and~$e\not=1$.

\mgni{\it Case 4.} Assume~$i=s\b+e$
and~$j=t\b+f$ with~$t<s-1$ and~$e,f\not=1$. Then we have
$$\eqalign{\xx&_{s\b+e}(\xx_1\cdots\xx_p)^{(t+2)\a}
=^{{}^{\!\!\!\!\!\!\!\rm\scr(2)}}\!
(\xx_1\cdots\xx_p)^{(t+2)\a}\xx_{(s-t-2)\b+e}\cr
&=^{{}^{\!\!\!\!\!\!\!\rm\scr(1)}}\!
(\xx_1\cdots\xx_p)^{t\a}
\xx_f\xx_{f+1}\cdots\xx_p(\xx_1\cdots\xx_p)^\a
\xx_{1+p-\b}\cdots\xx_{f-1+p-\b}
(\xx_1\cdots\xx_p)^{\a-1}\xx_{(s-t-2)\b+e}\cr
&=^{{}^{\!\!\!\!\!\!\!\rm\scr(2)}}\!
\xx_{t\b+f}(\xx_1\cdots\xx_p)^{t\a}
\xx_{f+1}\cdots\xx_p(\xx_1\cdots\xx_p)^\a
\xx_{1+p-\b}\cdots\xx_{f-1+p-\b}
(\xx_1\cdots\xx_p)^{\a-1}\xx_{(s-t-2)\b+e}.}$$
This provides the expression
of~$\xx_{s\b+e}\dR\xx_{t\b+f}$ for either~$0<t+1<s$
or~$0<s+1<t$ and~$e,f\not=1$.

\mgni{\it Case 5.} Assume~$i=s\b+e$
and~$j=t\b+f$ with~$t=s-1$ and~$e,f\not=1$. Then we have
$$\eqalign{\xx_{s\b+e}&(\xx_1\cdots\xx_p)^{s\a}\xx_{e+1}\cdots\xx_p
(\xx_1\cdots\xx_p)^\a\xx_{1+p-\b}\cdots\xx_{e-1+p-\b}\cr
&=^{{}^{\!\!\!\!\!\!\!\rm\scr(2)}}\!
(\xx_1\cdots\xx_p)^{s\a}\xx_e\xx_{e+1}\cdots\xx_p
(\xx_1\cdots\xx_p)^\a\xx_{1+p-\b}\cdots\xx_{e-1+p-\b}\cr
&=^{{}^{\!\!\!\!\!\!\!\rm\scr(1)}}\!
(\xx_1\cdots\xx_p)^{(s+1)\a+1}=(\xx_1\cdots\xx_p)^{(t+2)\a+1}\cr
&=^{{}^{\!\!\!\!\!\!\!\rm\scr(1)}}\!
(\xx_1\cdots\xx_p)^{t\a}\xx_f\xx_{f+1}\cdots\xx_p
(\xx_1\cdots\xx_p)^\a\xx_{1+p-\b}\cdots\xx_{f-1+p-\b}\cr
&=^{{}^{\!\!\!\!\!\!\!\rm\scr(2)}}\!
\xx_{t\b+f}(\xx_1\cdots\xx_p)^{t\a}\xx_{f+1}\cdots\xx_p
(\xx_1\cdots\xx_p)^\a\xx_{1+p-\b}\cdots\xx_{f-1+p-\b}.}$$
This provides the expression
of~$\xx_{s\b+e}\dR\xx_{t\b+f}$ for either~$t=s-1$
or~$s=t-1$ and~$e,f\not=1$.

\mnni{\it Case 6.} Assume~$i=s\b+e$
and~$j=s\b+f$ with~$1\not=e\not=f\not=1$. Then the computations are
similar with the previous ones. This completes the
proof.\EndPr 

\Sec Torus link groups are Garside groups

\bnni The main result of this section is that all torus link
groups are groups of fractions of Garside monoids,
which, as a consequence, provides explicit biautomatic
structures for each of them.

\Prop\Ç\thinGaussÈ{\sl
Let~$2\leq p\leq q$. The monoid~$\atl{p,q}$ is a Garside monoid
with Garside element~$\displaystyle\D=
(\xx_1\cdots\xx_p)^{q\over{p\sswedge q}}$.}

\EPr For~$p\geq 2$ and~$k\geq 1$, $\atl{p,pk}$ is
a Garside monoid (see~[\\\Dfx, Example~5]
\hbox{[\\\Pic, Chapter~VI])} with Garside element~$(\xx_1\cdots\xx_p)^k$.

From now on, we assume~$2\leq p\leq q$
and~$p\notdiv q$. We shall use the Garsidity
criterion of~[\\\Dhh] described in Section~\SecBack.
First, Adjan's criterion [\\\Adj] (see
also~[\\\Rem][\\\Hig]) allows us to claim that,
for~$2\leq p\leq q$ and~$p\notdiv q$, the
monoid~$\atl{p,q}$ is cancellative~; indeed, the
graphs~$LG$Ç\ArtinTorIIÈ{ }and~$RG$Ç\ArtinTorIIÈ{ }are star graphs
centered respectively in~$\xx_1$ and in~$\xx_p$, hence cycle free. We
define~$\a,\b$ to be respectively the integer part and the remainder
of~$q\over p$. In particular, we have~$q=p\a+\b$ with~$1\leq\b\leq p-1$.

\mgni Next, in order to show that~$\displaystyle\D=
(\xx_1\cdots\xx_p)^{q\over{p\sswedge q}}$ is a
Garside element for~$\atl{p,q}$, we consider the
relations of the presentation~Ç\originalÈ
for~$\atl{p,q}$, which can be rewritten as
$$\eqalignno{\xx_i(\xx_1\cdots\xx_p)^{\a+1}
&=(\xx_1\cdots\xx_p)^{\a+1}\xx_{i+p-\b}
&\hbox{for~$1\leq i\leq\b$,\quad}\cr
\xx_i(\xx_1\cdots\xx_p)^\a
&=(\xx_1\cdots\xx_p)^\a\xx_{i-\b}
&\hbox{for~$\b+1\leq i\leq p$.\quad}}$$Let~$[z]$
denote the integer part of~$z$. The
map~$\tau$ of~$\{1,\ldots,\xx_p\}$ defined
by$$\tau(i)=p[{{q-i}\over p}]+p-q+i
=\cases{i+p-\b&for~$1\leq i\leq\b$,\cr
i-\b&for~$\b+1\leq i\leq p$,}$$ is a permutation
of~{\go S}$_p$ with order~${p\over{p\swedge q}}$.
More precisely, let us consider the
braid~$\tau_{\lambda,\mu}$ defined as the positive
$\lambda\!+\!\mu$-strand braid where the~$\mu$
strands initially at positions~$\lambda+1$
to~$\lambda+\mu$ cross over the first~$\lambda$
strands~:

\midinsert
\centerline{
\newlabel <-13, 12> $\tau_{\lambda,\mu}$
\newlabel <1.1, 26.4> $\overbrace{\hskip16.2mm}$
\newlabel <20, 26.5> $\overbrace{\hskip26.3mm}$
\newlabel <8.2, 29.5> $\lambda$
\newlabel <32, 29.5> $\mu$
\epsfscale=500
\epsfbox{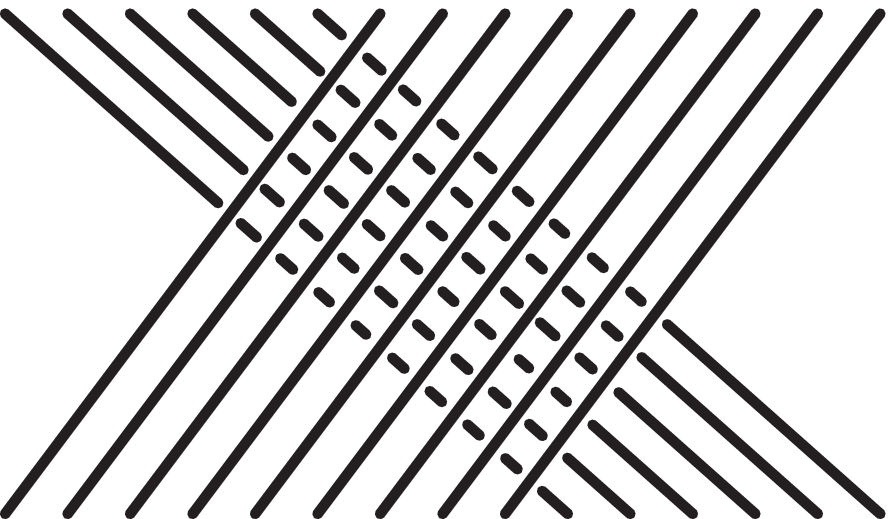}
}
\endinsert

The minimal integer~$k>0$
for~$\tau_{\lambda,\mu}^k$ to be a pure braid
is~$k={{\lambda+\mu}\over{\lambda\swedge\mu}}$, and,
moreover, in the pure braid~$\tau_{\lambda,\mu}^k$,
every strand crosses over (\resp crosses under) the
others~$\mu\over{\lambda\swedge\mu}$ times (\resp
$\lambda\over{\lambda\swedge\mu}$ times). As~$\tau$ is the
permutation associated with~$\tau_{\b,p-\b}$, we deduce
that the minimal integer~$h>0$ such
that~$\xx_i(\xx_1\cdots\xx_p)^h=(\xx_1\cdots\xx_p)^h\xx_i$
holds for every~$1\leq i\leq p$
is$$h={\b\over{\b\swedge(p-\b)}}(\a+1)
+{{p-\b}\over{\b\swedge(p-\b)}}\a={q\over{p\swedge
q}},$$ which proves that~$\displaystyle\D=
(\xx_1\cdots\xx_p)^{q\over{p\sswedge q}}$ is a
Garside element for~$\atl{p,q}$.

\mgni Finally, the presentation of Proposition~Ç\ComplementedÈ is
normed, as the relations preserve the length. So, we are left
with checking the local coherence, \ie, with checking
that
$$((\xx_i\dR\xx_j)\dR(\xx_i\dR\xx_k))
\dR((\xx_j\dR\xx_i)\dR(\xx_j\dR\xx_k))=\e$$
holds for~$1\leq i\not=j\not=k\leq p$. The verification
is easy but somehow tedious~: the
decomposition into some 133 subcases and the systematic
verification are given in Appendix.\EndPr

\bgni Gathering Propositions~Ç\AutomaticÈ
and~Ç\thinGaussÈ, we obtain~:

\ECor\Ç\TLGAutomaticÈ{\sl Every torus
link group admits an explicit biautomatic structure, given by the
associated torus link Artin monoid.}

\EEx{ Figure~{\FigWirtFourSix} displays the lattice of simples
in the monoid~$\atl{4,6}$.}

\midinsert
\centerline{
\newlabel <53.7, -3.4> $1$
\newlabel <53.2, 126> $\D$
\newlabel <31, 9> $\xx_2$
\newlabel <43, 11> $\xx_1$
\newlabel <63, 11> $\xx_3$
\newlabel <75.2, 9> $\xx_4$
\epsfscale=565
\epsfbox{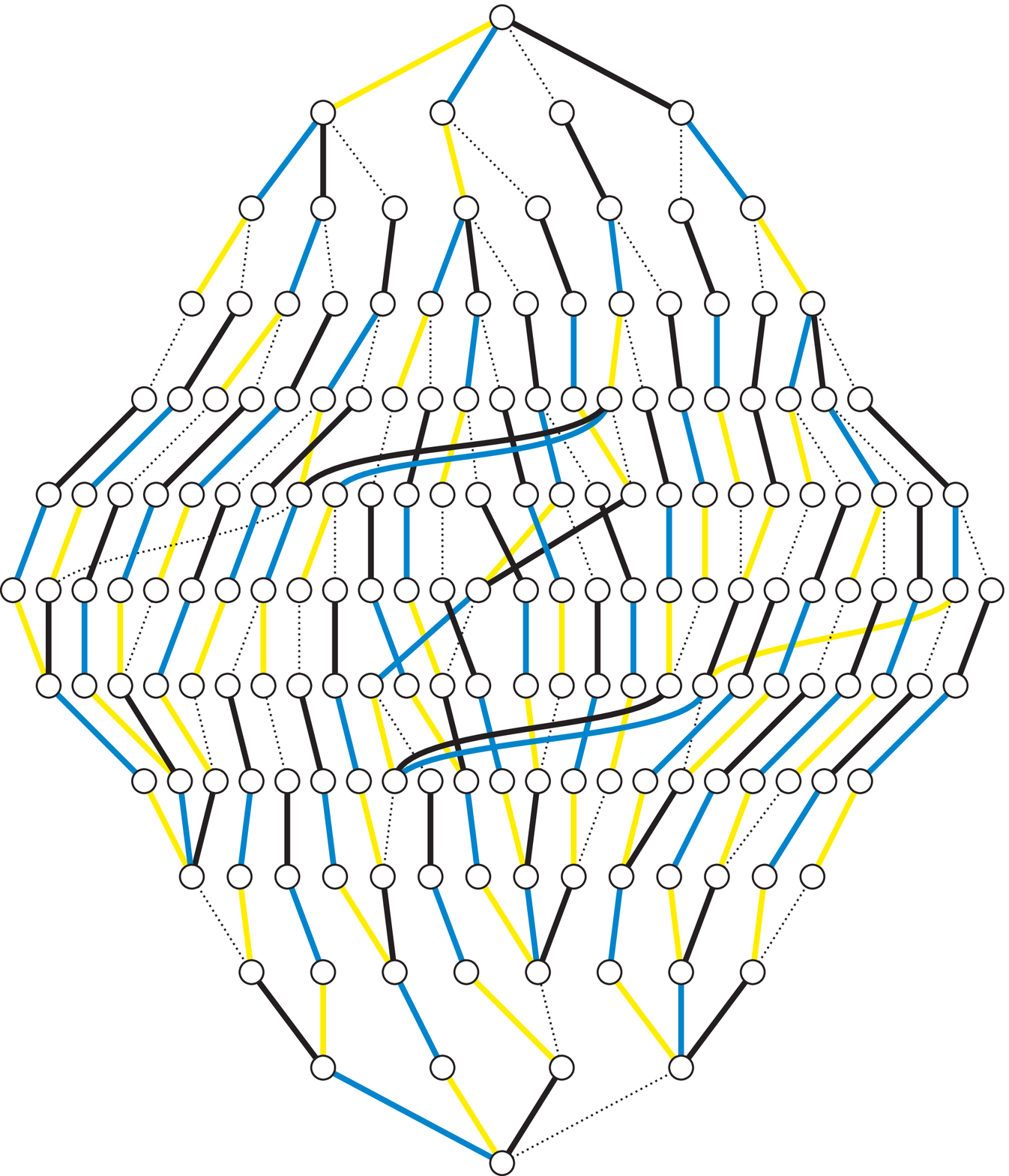}
}
\vglue 3mm
\centerline{\IX{\bfIX
Figure~\FigWirtFourSix.} The lattice of simples
in~$\atl{4,6}$.}
\endinsert

\ERem{ A result of Burde-Zieschang [\\\BuZ,
theorem~6.1] claims that a non-trivial knot whose group has a
non-trivial center is a torus knot. By comparing this
result with~[\\\Pib, Proposition~4.1], which claims that every Garside groups is an iterated crossed product of Garside groups
with a cyclic center generated by some power of its element~$\D$, we
deduce that, among all knots, only the torus knots---and the trivial
knot---have a group which is a Garside
group.}

\Sec Other  monoids
for torus link groups~?

\bnni Like every Garside group, a torus link group can
be the group of fractions of several (Garside or not)
monoids. For instance, we know that, for~$p,q>1$ coprime, the
$(p,q)$-torus knot group admits the (monoid)
presentation~$\langle~\xx,\yy:\xx^p=\yy^q~\rangle$, and, as
already mentionned in~[\\\Dfx], the monoid admitting this
presentation is a Garside monoid. In this final
section, we consider other monoid presentations for
torus link groups, inspired by the Wirtinger group
presentations of a link group.

\font\go=eufm10

\def\MoinsUn{\;\!\!}

\EDef\Ç\WirtMonoidÈ{ Let~$\pi$ be a regular projection of an
oriented link. We call \hbox{\it Wirtinger $\pi$-monoid}
the monoid~$\wir{\pi}$ admitting the
presentation $\langle~X:R~\rangle$, where~$X$
is the set of Wirtinger generators, \ie, the set
of all overcrossing arcs of~$\pi$, and $R$ is the set of
positive relations associated with the crossings
following the scheme of Figure~\FigRelation.}

\midinsert
\centerline{
\newlabel <21.7, 21> $\xx_i$
\newlabel <7, 9.5> $\xx_j$
\newlabel <29, 9.5> $\xx_k$
\epsfscale=500
\epsfbox{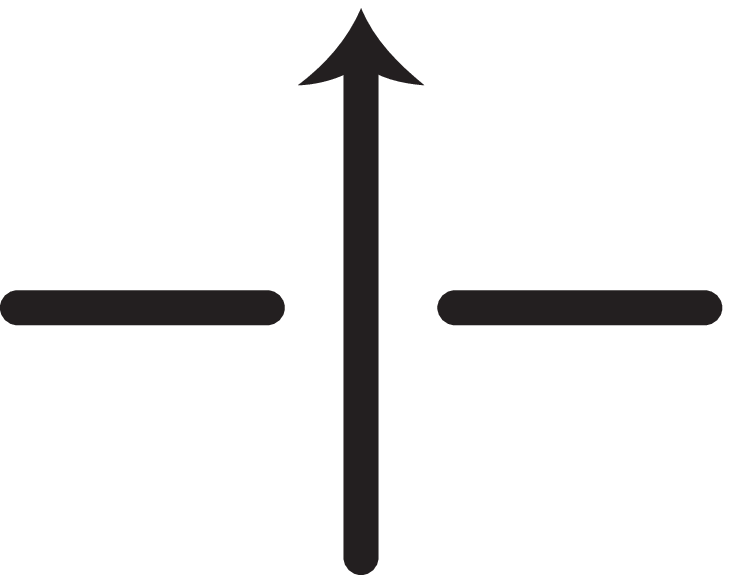}
}
\vglue 3mm
\centerline{\IX{\bfIX Figure~\FigRelation.} To this crossing corresponds
the relation~$\xx_i\xx_j=\xx_k\xx_i$.}
\endinsert

\ni Let~$\wir{p,q}$ denote the Wirtinger monoid
associated with the diagram
of the closure of the braid
word~$(\s_1\ldots\s_{p-1})^q$, where the~$\s_i$'s
are the Artin generators
(see~Section~{\SecPres}).

\bgni It turns out that certain projections~$\pi$ of torus
links yield a Wirtinger $\pi$-monoid that is a Garside monoid, and, therefore, yield an explicit
automatic structure for the associated torus link
group. The point is that the automatic structure so-obtained
can be simpler that the general one described in the
previous section.

\mgni For every~$q>1$ odd, the monoid~$\wir{2,q}$
is a Garside monoid, whose group of fractions is
isomorphic to the type~$I_2(q)$  Artin group. In particular,
the group of the
$(2,3)$-torus knot---\ie, the trefoil knot---is the braid
group~$B_3$.

In contradistinction, the
monoids~$\wir{3,2}$ and~$\wir{4,2}$ are not Garside. The
monoid~$\wir{3,2}$ admits the presentation$$\eqalign{
\langle~\s_1,\s_2,\tau_1,\tau_2:
~&\tau_2\s_1=\s_1\s_2=\s_2\tau_2,~
\s_1\tau_1=\s_2\s_1=\tau_1\s_2~\rangle,}$$ 
and, in particular, the
elements~$\s_1$ and~$\s_2$ do not admit a
unique lcm (neither on the right, nor on the left).
Nevertheless, $\wir{3,2}$ admits an
element~$\Delta$ whose left and right divisors
coincide and generate~$\wir{3,2}$,
namely~$\Delta=\tau_2^{}\s_1^2$. This implies
that any two elements in~$\wir{3,2}$ admit
a right (and a left) common multiple.
Moreover, adding the
relations~$\tau_1^{}\s_2^2
=\tau_2^{}\s_1^2=\s_1^2\tau_1^{}=\s_2^2\tau_2^{}$
to the previous presentation allows to obtain a {\it complete}
presentation for the monoid~$\wir{3,2}$ which, according to a
criterion of~[\\\Djj], implies that~$\wir{3,2}$ is
cancellative. In particular,
$\wir{3,2}$ satisfies Ore's conditions, hence embeds into its
group of fractions, namely~$B_3$. Let us mention that this
embedding result answers positively to a question of Han
\&~Ko in~[\\\HKo] about this Sergiescu monoid  for~$B_3$.
Figure~{\FigWirtThreeTwo} displays the Cayley graph
of~$\wir{3,2}$ restricted to divisors of~$\Delta$.

The monoid~$\wir{4,2}$ does not embed into a group, for,
in particular, it is no cancellative~: indeed, $\wir{4,2}$ is
$$\eqalign{\langle~\s_1,\s_2,\tau_1,\tau_2,\rho_1,\rho_2:
&~\s_1\s_2=\s_2\rho_1=\tau_2\s_1,~\s_1\tau_1=\rho_1\s_1,\cr
&~\s_2\s_1=\s_1\rho_2=\tau_1\s_2,~\s_2\tau_2=\rho_2\s_2~\rangle,}$$
and both~$\s_1\s_2\s_1=\s_1\rho_2\tau_1$
and~$\s_2\s_1\not=\rho_2\tau_1$ hold in~$\wir{4,2}$.

\midinsert
\centerline{
\newlabel <25.7, -3.4> $1$
\newlabel <-2, 13> $\tau_{\MoinsUn\scr2}$
\newlabel <14, 13> $\s_{\MoinsUn\scr1}$
\newlabel <37, 13> $\s_{\MoinsUn\scr2}$
\newlabel <53, 13> $\tau_{\MoinsUn\scr1}$
\newlabel <25.2, 42.3> $\Delta$
\epsfscale=500
\epsfbox{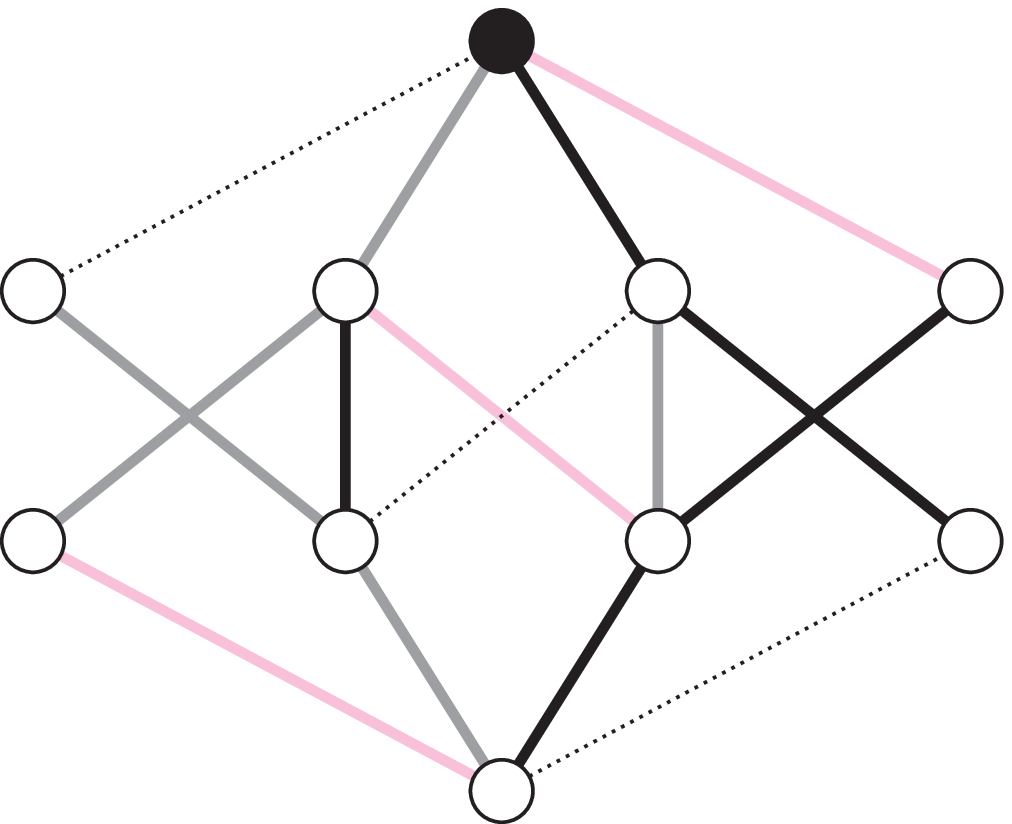}
}
\vglue 3mm
\centerline{\IX{\bfIX Figure~\FigWirtThreeTwo.} The 
monoid~$\wir{3,2}$ is not Garside, but it embeds into its
group of fractions~$B_3$.}
\endinsert

\mgni The monoid~$\wir{3,4}$ is a Garside monoid; indeed, it is
the monoid
$$\eqalign{\langle~\xx_1,\ldots,\xx_8:~&\xx_1\xx_7=\xx_7\xx_2=\xx_6\xx_1,~
\xx_2\xx_5=\xx_5\xx_3=\xx_3\xx_6,\cr
&\xx_1\xx_4=\xx_3\xx_1=\xx_8\xx_3,~
\xx_4\xx_7=\xx_5\xx_8=\xx_7\xx_5
~\rangle,}$$which admits the normed coherent complemented
presentation
$$\eqalign{
\langle~\xx_1,&\ldots,\xx_8:
\xx_1\xx_4\xx_7=\xx_2\xx_5\xx_1,~
\xx_1\xx_4=\xx_3\xx_1,~
\xx_1\xx_4\xx_7=\xx_4\xx_7\xx_3,\cr
&\xx_1\xx_4\xx_7=\xx_5\xx_1\xx_4,~   
\xx_1\xx_7=\xx_6\xx_1,~   
\xx_1\xx_7=\xx_7\xx_2,~
\xx_1\xx_4=\xx_8\xx_3,~
\xx_2\xx_5=\xx_3\xx_6,\cr
&\xx_2\xx_5\xx_1=\xx_4\xx_7\xx_3,~
\xx_2\xx_5=\xx_5\xx_3,~
\xx_2\xx_5\xx_1=\xx_6\xx_1\xx_5,~
\xx_2\xx_5\xx_1=\xx_7\xx_2\xx_5,\cr
&\xx_2\xx_5\xx_1=\xx_8\xx_3\xx_7,~
\xx_3\xx_1\xx_7=\xx_4\xx_7\xx_3,~
\xx_3\xx_6=\xx_5\xx_3,~
\xx_3\xx_1\xx_7=\xx_6\xx_1\xx_5,\cr
&\xx_3\xx_1\xx_7=\xx_7\xx_2\xx_5,~
\xx_3\xx_1=\xx_8\xx_3,~
\xx_4\xx_7=\xx_5\xx_8,~
\xx_4\xx_7\xx_3=\xx_6\xx_1\xx_5,\cr
&\xx_4\xx_7=\xx_7\xx_5,~
\xx_4\xx_7\xx_3=\xx_8\xx_3\xx_7,~
\xx_5\xx_1\xx_4=\xx_6\xx_1\xx_5,~
\xx_5\xx_8=\xx_7\xx_5,\cr
&\xx_5\xx_1\xx_4=\xx_8\xx_3\xx_7,~
\xx_6\xx_1=\xx_7\xx_2,~
\xx_6\xx_1\xx_5=\xx_8\xx_3\xx_7,~
\xx_7\xx_2\xx_5=\xx_8\xx_3\xx_7~\rangle.}$$
A Garside element is $\D=\xx_1\xx_4\xx_7$,
and~$\wir{3,4}$ is a Garside monoid, whose center is
generated by the element~$\D^4$. The
lattice of simple elements in~$\wir{3,4}$ is displayed
on~Figure~\FigWirtThreeFour. In contradistinction, the Wirtinger
monoid~$\wir{4,3}$ is not Garside.

\midinsert
\centerline{
\newlabel <54.9, -3.4> $1$
\newlabel <-2.5, 16> $\xx_{\MoinsUn\scr1}$
\newlabel <12.5, 16> $\xx_{\MoinsUn\scr2}$
\newlabel <27.8, 16> $\xx_{\MoinsUn\scr3}$
\newlabel <43, 16> $\xx_{\MoinsUn\scr4}$
\newlabel <65.7, 16> $\xx_{\MoinsUn\scr5}$
\newlabel <81, 16> $\xx_{\MoinsUn\scr6}$
\newlabel <96.3, 16> $\xx_{\MoinsUn\scr7}$
\newlabel <111.6, 16> $\xx_{\MoinsUn\scr8}$
\newlabel <54.5, 72.1> $\D$
\epsfscale=600
\epsfbox{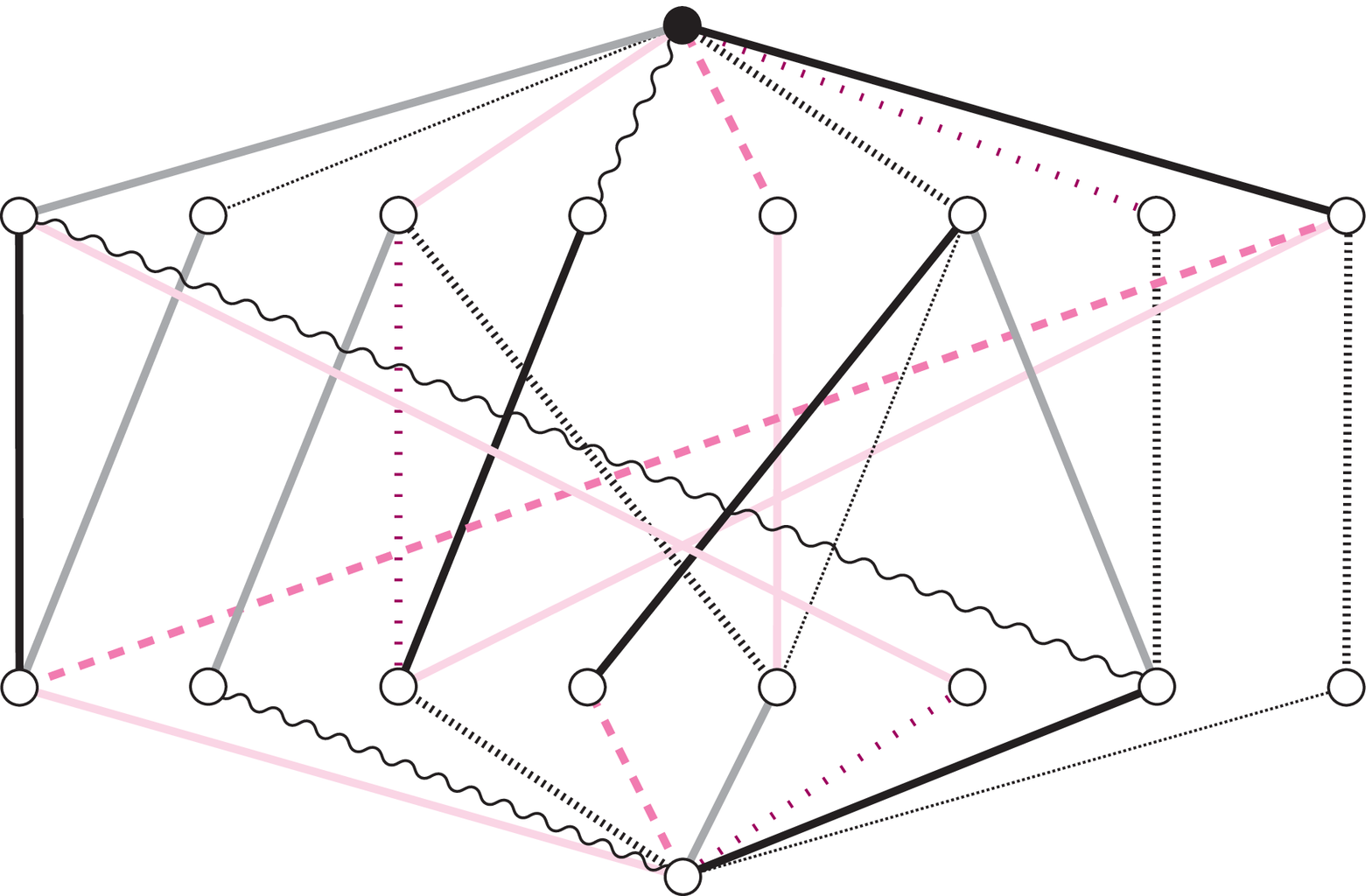}
}
\vglue 3mm
\centerline{\IX{\bfIX
Figure~\FigWirtThreeFour.} The lattice of simples
in~$\wir{3,4}$.}
\endinsert

\ERem{ There are many other examples of Garside monoids
whose groups of fractions are isomorphic to a torus link
group. For instance, the group of fractions of the
monoid~$\Knuth$ of~Example~Ç\ExKnuthÈ is isomorphic to the
$(3,4)$-torus knot group. For every~$k\geq 0$, the
monoid~$M^{\scr\!(2,2)}_k$ presented
by$$\langle~\xx,\yy:\xx\yy^2\xx\cdot
(\yy\xx)^k\yy\cdot\xx\yy^2\xx=
(\yy\xx)^k\yy\cdot\xx\yy^2\xx\cdot (\yy\xx)^k\yy~\rangle$$is
a Garside monoid, whose Garside element is the square
of the lcm of the atoms. The lattice of simples
in~$M^{\scr\!(2,2)}_0$ is displayed in Figure~\FigTypeII. The
groups of fractions of the monoids~$M^{\scr\!(2,2)}_k$ are
all isomorphic to the (3,4)-torus knot
group~: take~$\zz=\xx\yy^2(\xx\yy)^{k+1}$
and~$\ttt=\yy(\xx\yy)^{k+1}$,
\ie, $\xx=(\zz\ttt^{-1})^{k+2}\ttt^{-1}$
and~$\yy=\ttt(\ttt\zz^{-1})^{k+1}$.}

\midinsert
\centerline{
\newlabel < 34.7, -3> $1$
\newlabel < 34.3, 83> $\D$
\epsfscale=350
\epsfbox{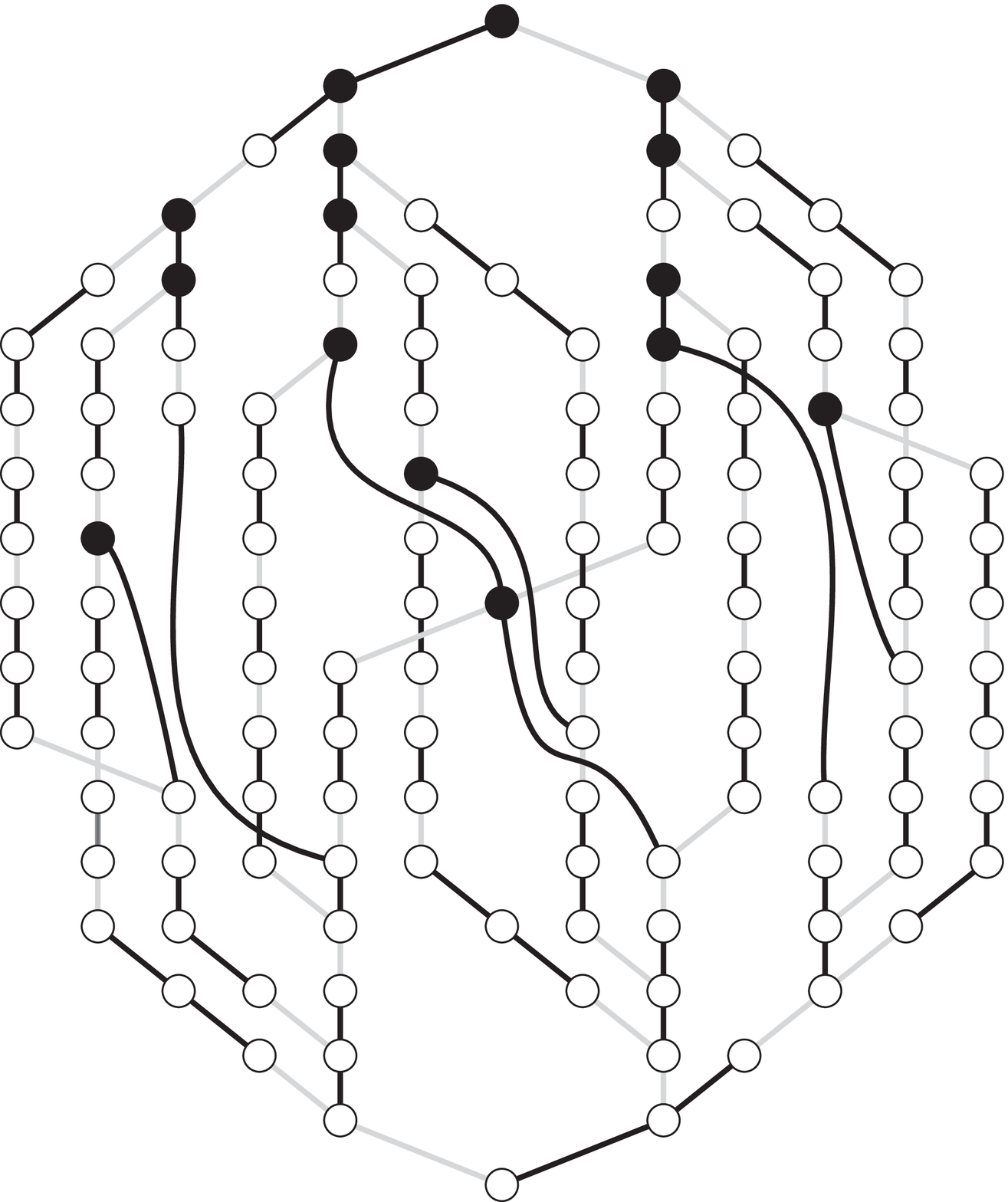}
}
\vglue 3mm
\centerline{\IX{\bfIX Figure~\FigTypeII.} The
lattice of simples in~$\langle~\xx,\yy:
\xx\yy\yy\xx\yy\xx\yy\yy\xx
=\yy\xx\yy\yy\xx\yy~\rangle$.}
\endinsert

\mgni Finally, $\wir{4,6}$ is a Garside monoid~: the
lattice of the 68~simples in~$\wir{4,6}$ is displayed on
Figure~\FigWirtFourSix, and can be compared with the lattice
of the 176~simples in~$\atl{4,6}$ on~Figure~\FigArtFourSix.

\midinsert
\vglue 0.7cm
\centerline{
\newlabel <56.8, -3.4> $1$
\newlabel <56.2, 67> $\D$
\epsfscale=400
\epsfbox{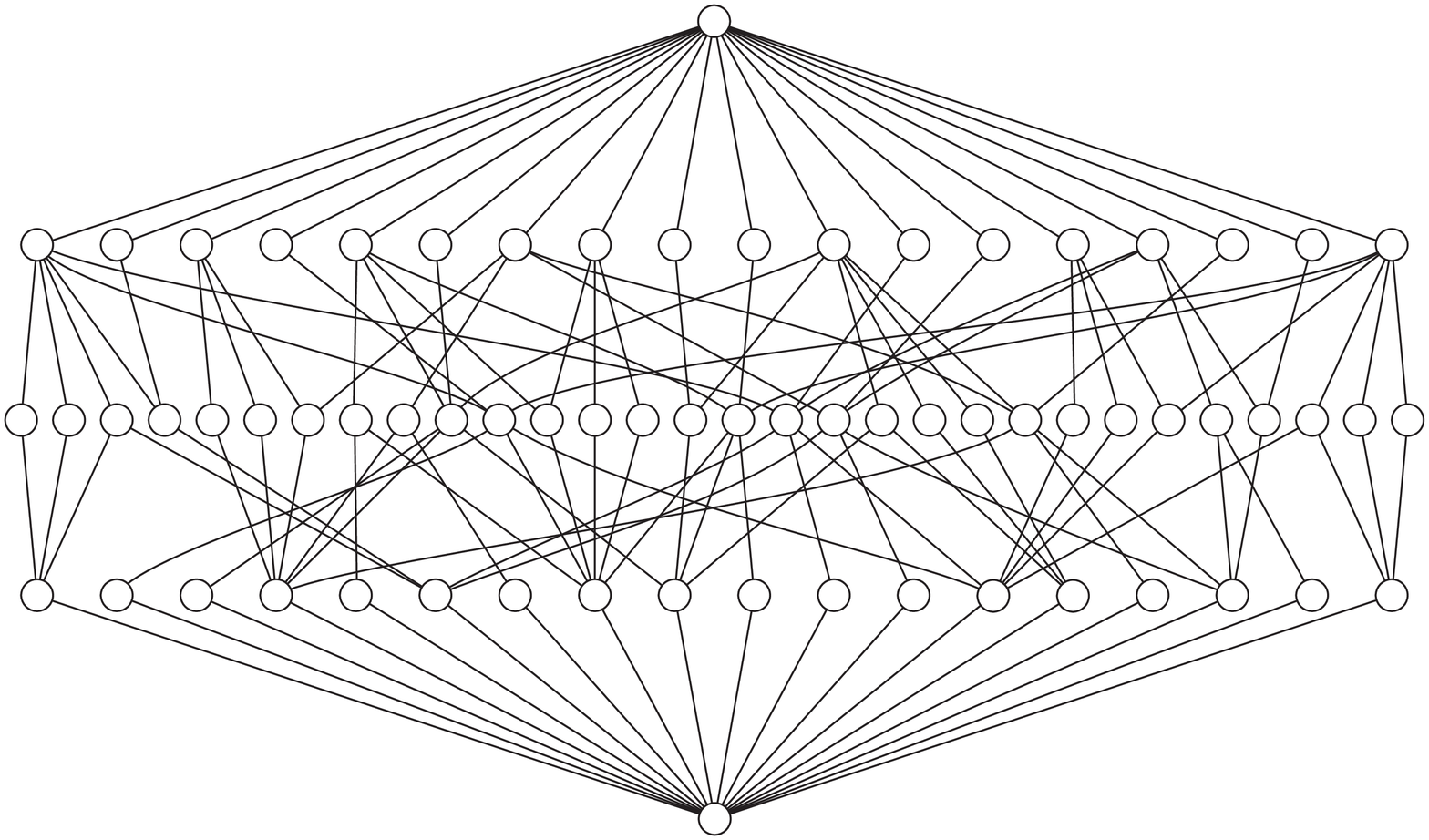}
}
\vglue 3mm
\centerline{\IX{\bfIX
Figure~\FigWirtFourSix.} The (Hass diagram of the) lattice of
simples in~$\wir{4,6}$.}
\vglue 0.7cm
\endinsert

\bgni We thus are naturally led to the two following
problems~:

\mgni {\bf Question 1. }{\sl Among all possible Wirtinger
monoids of torus links, which ones are Garside?}

\mgni {\bf Question 2. }{\sl Among all possible Wirtinger
monoids of links, which ones embed into a~group?}

\vfill\eject
\bg\bgni\centerline{\scXII
Appendix}

\bnni Here is the full detail of verification of 
local coherence of the presentation of
Proposition~Ç\ComplementedÈ, completing the proof of
Proposition~Ç\thinGaussÈ. The notations here are the ones used
in~Proposition~Ç\ComplementedÈ. Above all, let us observe
that, from the relation
$$(\xx_1\cdots\xx_p)\cdot
(\xx_1\cdots\xx_p)^{\a-1}\xx_{(z-1)\b+1}
=\xx_{z\b+1}\cdot(\xx_1\cdots\xx_p)^\a$$
of Presentation~Ç\ArtinTorIIÈ, we deduce
$$\eqalign{w_{\scr000}\dR\xx_{z\b+1}&=w_{\scr00(\a-2)}\xx_{(z-1)\b+1}\cr
\xx_{z\b+1}\dR w_{\scr000}&=w_{\scr00(\a-1)}.}$$
Let~$i=s\b+e$, $j=t\b+f$, $k=z\b+g$
and~$u_{ijk}=((\xx_i\dR\xx_j)\dR(\xx_i\dR\xx_k))
\dR((\xx_j\dR\xx_i)\dR(\xx_j\dR\xx_k))$. By~definition
of~$\dR$, coherence always holds for triplets~$(u,v,w)$ where
two among~$u,v,w$ are equal~: in the
following, we shall assume
$s=t\Rightarrow e\not=f$,
$t=z\Rightarrow f\not=g$ and~$z=s\Rightarrow g\not=e$.

\def\CaseI{1}
\def\CaseII{2}
\def\CaseIII{3}
\def\CaseIV{4}
\def\CaseV{5}
\def\CaseVI{6}
\def\CaseVII{7}
\def\CaseVIII{8}

\mgni{\it Case~\CaseI:} $e=1$, $f=1$, $g=1$.

\ni{\it Subcase~\CaseI.1:} $z<t<s$.
We find$$\eqalign{u_{ijk}&=(w_{\scr(t\a)0(\a-1)}\dR
w_{\scr(z\a)0(\a-1)})
\dR(w_{\scr(t\a)1(\a-1)}\xx_{(s-t-1)\b+1}
\dR w_{\scr(z\a)0(\a-1)})=\e\dR\e=\e.}$$
{\it Subcase~\CaseI.2:} $t<z<s$.
We find$$\eqalign{u_{ijk}&=(w_{\scr(t\a)0(\a-1)}\dR
w_{\scr(z\a)0(\a-1)})
\dR(w_{\scr(t\a)1(\a-1)}\xx_{(s-t-1)\b+1}
\dR
w_{\scr(t\a)1(\a-1)}\xx_{(z-t-1)\b+1})\cr
&=(\xx_1\cdots\xx_p)^{(z-t)\a}
\dR(\xx_{(s-t-1)\b+1}
\dR\xx_{(z-t-1)\b+1})\cr
&=(\xx_1\cdots\xx_p)^{(z-t)\a}
\dR w_{\scr((z-t-1)\a)0(\a-1)}=\e.}$$
{\it Subcase~\CaseI.3:} $t<s<z$.
We find$$\eqalign{u_{ijk}&=(w_{\scr(t\a)0(\a-1)}\dR
w_{\scr(s\a)1(\a-1)}\xx_{(z-s-1)\b+1})\cr
&\phantom{=5}\dR(w_{\scr(t\a)1(\a-1)}\xx_{(s-t-1)\b+1}
\dR w_{\scr(t\a)1(\a-1)}\xx_{(z-t-1)\b+1})\cr
&=(w_{\scr((s-t-1)\a)1(\a-1)}\xx_{(z-s-1)\b+1})\dR(\xx_{(s-t-1)\b+1}
\dR\xx_{(z-t-1)\b+1})=\e.}$$
{\it Subcase~\CaseI.4:} $z<s<t$.
We find$$\eqalign{u_{ijk}&=(w_{\scr(s\a)1(\a-1)}\xx_{(t-s-1)\b+1}\dR
w_{\scr(z\a)0(\a-1)})
\dR(w_{\scr(s\a)0(\a-1)}
\dR w_{\scr(z\a)0(\a-1)})=\e\dR\e=\e.}$$
{\it Subcase~\CaseI.5:} $s<z<t$.
We find$$\eqalign{u_{ijk}&=(w_{\scr(s\a)1(\a-1)}\xx_{(t-s-1)\b+1}\dR
w_{\scr(s\a)1(\a-1)}\xx_{(z-s-1)\b+1})
\dR(w_{\scr(s\a)0(\a-1)}
\dR w_{\scr(z\a)0(\a-1)})\cr
&=(\xx_{(t-s-1)\b+1}\dR\xx_{(z-s-1)\b+1})\dR
w_{\scr((z-s-1)\a)0(\a-1)}=\e.}$$
{\it Subcase~\CaseI.6:} $s<t<z$.
We
find$$\eqalign{u_{ijk}&=(w_{\scr(s\a)1(\a-1)}\xx_{(t-s-1)\b+1}\dR
w_{\scr(s\a)1(\a-1)}\xx_{(z-s-1)\b+1})
\dR(w_{\scr(s\a)0(\a-1)}\cr
&\phantom{=5}\dR
w_{\scr(t\a)1(\a-1)}\xx_{(z-t-1)\b+1})\cr
&=(\xx_{(t-s-1)\b+1}\dR\xx_{(z-s-1)\b+1})
\dR(w_{\scr((t-s-1)\a)1(\a-1)}\xx_{(z-t-1)\b+1})=\e.}$$

\mgni{\it Case~\CaseII:} $e\not=1$, $f=1$,
$g=1$.

\ni{\it Subcase~\CaseII.1:} $z<t<s$.
We find$$\eqalign{u_{ijk}&=(w_{\scr(t\a)0(\a-1)}\dR
w_{\scr(z\a)0(\a-1)})
\dR(w_{\scr(t\a)1(\a-1)}\xx_{(s-t-1)\b+e}
\dR w_{\scr(z\a)0(\a-1)})=\e\dR\e=\e.}$$
{\it Subcase~\CaseII.2:} $t<z<s$.
We find$$\eqalign{u_{ijk}
&=(w_{\scr(t\a)0(\a-1)}\dR w_{\scr(z\a)0(\a-1)})
\dR(w_{\scr(t\a)1(\a-1)}\xx_{(s-t-1)\b+e}
\dR w_{\scr(t\a)1(\a-1)}\xx_{(z-t-1)\b+1})\cr
&=w_{\scr((z-t-1)\a)0(\a-1)}\dR(\xx_{(s-t-1)\b+e}\dR\xx_{(z-t-1)\b+1})=\e.}$$
{\it Subcase~\CaseII.3:} $t<s\leq z\leq s+1$.
We find$$\eqalign{u_{ijk}
&=(w_{\scr(t\a)0(\a-1)}\dR w_{\scr(s\a)e\a})
\dR(w_{\scr(t\a)1(\a-1)}\xx_{(s-t-1)\b+e}
\dR w_{\scr(t\a)1(\a-1)}\xx_{(z-t-1)\b+1})\cr
&=w_{\scr((s-t-1)\a)e\a}\dR(\xx_{(s-t-1)\b+e}\dR\xx_{(z-t-1)\b+1})=\e.}$$
{\it Subcase~\CaseII.4:} $t<s<z-1$.
We find$$\eqalign{u_{ijk}
&=(w_{\scr(t\a)0(\a-1)}\dR
w_{\scr(s\a)e\a}w_{\scr00(\a-2)}\xx_{(z-s-2)\b+1})\cr
&\phantom{=5}\dR(w_{\scr(t\a)1(\a-1)}\xx_{(s-t-1)\b+e}
\dR w_{\scr(t\a)1(\a-1)}\xx_{(z-t-1)\b+1})\cr
&=w_{\scr((s-t-1)\a)e\a} w_{\scr00(\a-2)}\xx_{(z-s-2)\b+1}\dR
(\xx_{(z-s-2)\b+1})\dR\xx_{(s-t-1)\b+e})=\e.}$$
{\it Subcase~\CaseII.5:} $z<s\leq t\leq s+1$.
We find$$\eqalign{u_{ijk}&=(w_{\scr(s\a)e\a}\dR
w_{\scr(z\a)0(\a-1)})\dR
(w_{\scr(t\a)1(s\a-t\a+\a)}\dR
w_{\scr(z\a)0(\a-1)})=\e\dR\e=\e.}$$
{\it Subcase~\CaseII.6:} $t=s+1=z+1$.
We find$$\eqalign{u_{ijk}&=(w_{\scr(s\a)e\a}\dR
w_{\scr(s\a)e\a})\dR
(w_{\scr(t\a)1(s\a-t\a+\a)}\dR
w_{\scr(z\a)0(\a-1)})=\e\dR\e=\e.}$$
{\it Subcase~\CaseII.7:} $z=s+1=t+1$.
We find$$\eqalign{u_{ijk}&=(w_{\scr(s\a)e\a}\dR
w_{\scr(s\a)e\a})\dR
(w_{\scr(t\a)1(s\a-t\a+\a)}\dR
w_{\scr(t\a)1(\a-1)}\xx_{(z-t-1)\b+1})=\e\dR
(w_{\scr000}\dR\xx_1)=\e.}$$
{\it Subcase~\CaseII.8:} $s\leq t\leq s+1<z$.
We find$$\eqalign{u_{ijk}&=(w_{\scr(s\a)e\a}\dR
w_{\scr(s\a)e\a}w_{\scr00(\a-2)}\xx_{(z-s-2)\b+1})\dR
(w_{\scr(t\a)1(s\a-t\a+\a)}\dR
w_{\scr(t\a)1(\a-1)}\xx_{(z-t-1)\b+1})\cr
&=\cases{w_{\scr00(\a-2)}\xx_{(z-s-2)\b+1}
\dR(w_{\scr000}\dR\xx_{(z-t-1)\b+1})&for~$t=s$\cr
w_{\scr00(\a-2)}\xx_{(z-s-2)\b+1}\dR w_{\scr00(\a-2)}\xx_{(z-t-1)\b+1}&for~$t=s+1$}=\e.}$$
{\it Subcase~\CaseII.9:} $z+1<s+1<t$.
We find$$\eqalign{u_{ijk}&=(w_{\scr(s\a)e\a}w_{\scr00(\a-2)}
\xx_{(t-s-2)\b+1}\dR
w_{\scr(z\a)0(\a-1)})\dR(w_{\scr(s\a)0(2\a-1)}\dR
w_{\scr(z\a)0(\a-1)})=\e\dR\e=\e.}$$
{\it Subcase~\CaseII.10:} $s\leq z\leq s+1<t$.
We find$$\eqalign{u_{ijk}&=(w_{\scr(s\a)e\a}w_{\scr00(\a-2)}\xx_{(t-s-2)\b+1}\dR
w_{\scr(s\a)e\a})\dR(w_{\scr(s\a)0(2\a-1)}\dR
w_{\scr(z\a)0(\a-1)})=\e\dR\e=\e.}$$
{\it Subcase~\CaseII.11:} $s+1<z<t$.
We find$$\eqalign{u_{ijk}&=(w_{\scr(s\a)e\a}w_{\scr00(\a-2)}\xx_{(t-s-2)\b+1}\dR
w_{\scr(s\a)e\a}w_{\scr00(\a-2)}\xx_{(z-s-2)\b+1})\dR(w_{\scr(s\a)0(2\a-1)}\dR
w_{\scr(z\a)0(\a-1)})\cr
&=(\xx_{(t-s-2)\b+1}\dR
\xx_{(z-s-2)\b+1})\dR\e=\e.}$$
{\it Subcase~\CaseII.12:} $s+1<t<z$.
We find$$\eqalign{u_{ijk}&=(w_{\scr(s\a)e\a}w_{\scr00(\a-2)}\xx_{(t-s-2)\b+1}\dR
w_{\scr(s\a)e\a}w_{\scr00(\a-2)}\xx_{(z-s-2)\b+1})\cr
&\phantom{=5}\dR(w_{\scr(s\a)0(2\a-1)}\dR
w_{\scr(t\a)1(\a-1)}\xx_{(z-t-1)\b+1})\cr
&=(\xx_{(t-s-2)\b+1}\dR
\xx_{(z-s-2)\b+1})\dR w_{\scr((t-s-2)\a)1(\a-1)}\xx_{(z-t-1)\b+1}=\e.}$$

\mgni{\it Case~\CaseIII:} $e=1$,
$f\not=1$, $g=1$.

\ni{\it Subcase~\CaseIII.1:} $z<s<t$.
We find$$u_{ijk}=(w_{\scr(s\a)1(\a-1)}\xx_{(t-s-1)\b+f}
\dR w_{\scr(z\a)0(\a-1)}
)\dR(w_{\scr(s\a)0(\a-1)}
\dR w_{\scr(z\a)0(\a-1)})=\e\dR\e=\e.$$
{\it Subcase~\CaseIII.2:} $s<z<t$.
We find$$\eqalign{u_{ijk}&=(w_{\scr(s\a)1(\a-1)}\xx_{(t-s-1)\b+f}
\dR w_{\scr(s\a)1(\a-1)}\xx_{(z-s-1)\b+1}
)\dR(w_{\scr(s\a)0(\a-1)}
\dR w_{\scr(z\a)0(\a-1)})\cr
&=(\xx_{(t-s-1)\b+f}\dR\xx_{(z-s-1)\b+1}
)\dR w_{\scr((z-s)\a)0(\a-1)}=\e.}$$
{\it Subcase~\CaseIII.3:} $s<t\leq z\leq t+1$.
We find$$\eqalign{u_{ijk}&=(w_{\scr(s\a)1(\a-1)}\xx_{(t-s-1)\b+f}
\dR w_{\scr(s\a)1(\a-1)}\xx_{(z-s-1)\b+1}
)\dR(w_{\scr(s\a)0(\a-1)}
\dR w_{\scr(t\a)f\a})\cr
&=(\xx_{(t-s-1)\b+f}\dR\xx_{(z-s-1)\b+1})\dR w_{\scr((t-s-1)\a)f\a}=\e.}$$
{\it Subcase~\CaseIII.4:} $s<t<z-1$.
We find$$\eqalign{u_{ijk}&=(w_{\scr(s\a)1(\a-1)}\xx_{(t-s-1)\b+f}
\dR w_{\scr(s\a)1(\a-1)}\xx_{(z-s-1)\b+1})\cr
&\phantom{=5}\dR(w_{\scr(s\a)0(\a-1)}
\dR w_{\scr(t\a)f\a}w_{\scr00(\a-2)}\xx_{(z-t-2)\b+1})\cr
&=(\xx_{(t-s-1)\b+f}\dR\xx_{(z-s-1)\b+1})\dR w_{\scr((t-s-1)\a)f\a}w_{\scr00(\a-2)}\xx_{(z-t-2)\b+1}=\e.}$$
{\it Subcase~\CaseIII.5.} Assumer~$z<t\leq s\leq t+1$.
We find$$u_{ijk}=(w_{\scr(s\a)1(t\a-s\a+\a)}
\dR w_{\scr(z\a)0(\a-1)}
)\dR(w_{\scr(t\a)f\a}
\dR w_{\scr(z\a)0(\a-1)})=\e\dR\e=\e.$$
{\it Subcase~\CaseIII.6.} Assume either~$z=t+1=s+1$
or~$z=t=s-1$.
We find$$\eqalign{u_{ijk}&=(w_{\scr(s\a)1(t\a-s\a+\a)}
\dR w_{\scr(s\a)1(\a-1)}\xx_{(z-s-1)\b+1}
)\dR(w_{\scr(t\a)f\a}
\dR w_{\scr(t\a)f\a})\cr
&=(w_{\scr00((t-s)\a)}\dR\xx_{(z-s-1)\b+1})\dR\e=\e.}$$
{\it Subcase~\CaseIII.7:} $t\leq s\leq t+1<z$.
We find$$\eqalign{u_{ijk}&=(w_{\scr(s\a)1(t\a-s\a+\a)}
\dR w_{\scr(s\a)1(\a-1)}\xx_{(z-s-1)\b+1})\dR(w_{\scr(t\a)f\a}
\dR w_{\scr(t\a)f\a}w_{\scr00(\a-2)}\xx_{(z-t-2)\b+1})\cr
&=\cases{(\xx_1\cdots\xx_p\dR\xx_{(z-s-1)\b+1})\dR w_{\scr00(\a-2)}\xx_{(z-t-2)\b+1}& for~$t=s$\cr
(w_{\scr00(\a-2)}\dR\xx_{(z-s-1)\b+1})\dR
w_{\scr00(\a-2)}\xx_{(z-t-2)\b+1}&for~$t+1=s$}=\e.}$$
{\it Subcase~\CaseIII.8:} $z<t<s-1$.
We find$$u_{ijk}=(w_{\scr(t\a)0(2\a-1)}
\dR w_{\scr(z\a)0(\a-1)}
)\dR(w_{\scr(t\a)f\a}w_{\scr00(\a-2)}\xx_{(s-t-2)\b+1}
\dR w_{\scr(z\a)0(\a-1)})=\e\dR\e=\e.$$
{\it Subcase~\CaseIII.9:} $t\leq z\leq
t+1<s$.
We find$$u_{ijk}=(w_{\scr(t\a)0(2\a-1)}
\dR w_{\scr(z\a)0(\a-1)}
)\dR(w_{\scr(t\a)f\a}w_{\scr00(\a-2)}\xx_{(s-t-2)\b+1}
\dR w_{\scr(t\a)f\a})=\e\dR\e=\e.$$
{\it Subcase~\CaseIII.10:} $t+1<z<s$.
We find$$\eqalign{u_{ijk}&=(w_{\scr(t\a)0(2\a-1)}\dR
w_{\scr(z\a)0(\a-1)})\dR(w_{\scr(t\a)f\a}w_{\scr00(\a-2)}\xx_{(s-t-2)\b+1}
\dR w_{\scr(t\a)f\a}w_{\scr00(\a-2)}\xx_{(z-t-2)\b+1})\cr
&=w_{\scr((z-t-2)\a)0(\a-1)}\dR(\xx_{(s-t-2)\b+1}\dR\xx_{(z-t-2)\b+1})=\e.}$$
{\it Subcase~\CaseIII.11:} $t+1<s<z$.
We find$$\eqalign{u_{ijk}&=(w_{\scr(t\a)0(2\a-1)}
\dR w_{\scr(s\a)1(\a-1)}\xx_{(z-s-1)\b+1})\cr
&\phantom{=5}\dR(w_{\scr(t\a)f\a}w_{\scr00(\a-2)}\xx_{(s-t-2)\b+1}
\dR w_{\scr(t\a)f\a}w_{\scr00(\a-2)}\xx_{(z-t-2)\b+1})\cr
&=w_{\scr((s-t-2)\a)1(\a-1)}\xx_{(z-s-1)\b+1}\dR(\xx_{(s-t-2)\b+1}\dR\xx_{(z-t-2)\b+1})=\e.}$$

\mgni{\it Case~\CaseIV:} $e=1$,
$f=1$, $g\not=1$.

\ni{\it Subcase~\CaseIV.1:} $t<s<z$.
We find$$\eqalign{u_{ijk}&=(w_{\scr(t\a)0(\a-1)}
\dR w_{\scr(s\a)1(\a-1)}\xx_{(z-s-1)\b+g})\cr
&\phantom{=5}\dR(w_{\scr(t\a)1(\a-1)}\xx_{(s-t-1)\b+1}\dR w_{\scr(t\a)1(\a-1)}\xx_{(z-t-1)\b+g})\cr
&=(w_{\scr((s-t-1)\a)1(\a-1)}\xx_{(z-s-1)\b+g})\dR(\xx_{(s-t-1)\b+1}\dR\xx_{(z-t-1)\b+g})=\e.}$$
{\it Subcase~\CaseIV.2:} $t<z\leq s\leq z+1$.
We find$$\eqalign{u_{ijk}&=(w_{\scr(t\a)0(\a-1)}
\dR w_{\scr(s\a)1(z\a-s\a+\a)}
)\dR(w_{\scr(t\a)1(\a-1)}\xx_{(s-t-1)\b+1}
\dR w_{\scr(t\a)1(\a-1)}\xx_{(z-t-1)\b+g})\cr
&=w_{\scr((s-t-1)\a)1(z\a-s\a+\a)}\dR(\xx_{(s-t-1)\b+1}\dR\xx_{(z-t-1)\b+g})=\e.}$$
{\it Subcase~\CaseIV.3:} $t=z=s-1$.
We find$$\eqalign{u_{ijk}&=(w_{\scr(t\a)0(\a-1)}
\dR w_{\scr(s\a)1(z\a-s\a+\a)}
)\dR(w_{\scr(t\a)1(\a-1)}\xx_{(s-t-1)\b+1}
\dR w_{\scr(t\a)1(z\a-t\a+\a)})\cr
&=w_{\scr010}\dR(w_{\scr(t\a)1(\a-1)}\xx_1
\dR w_{\scr(t\a)1\a})=w_{\scr010}\dR(\xx_1
\dR\xx_1\cdots\xx_p)=\e.}$$
{\it Subcase~\CaseIV.4:} $t<z<s-1$.
We find$$\eqalign{u_{ijk}&=(w_{\scr(t\a)0(\a-1)}
\dR w_{\scr(z\a)0(2\a-1)}
)\dR(w_{\scr(t\a)1(\a-1)}\xx_{(s-t-1)\b+1}
\dR w_{\scr(t\a)1(\a-1)}\xx_{(z-t-1)\b+g})\cr
&=w_{\scr((z-t-1)\a)0(2\a-1)}\dR(\xx_{(s-t-1)\b+1}\dR\xx_{(z-t-1)\b+g})=\e.}$$
{\it Subcase~\CaseIV.5:} $z\leq t\leq z+1<s$.
We find$$\eqalign{u_{ijk}&=(w_{\scr(t\a)0(\a-1)}
\dR w_{\scr(z\a)0(2\a-1)}
)\dR(w_{\scr(t\a)1(\a-1)}\xx_{(s-t-1)\b+1}
\dR w_{\scr(t\a)1(z\a-t\a+\a)})\cr
&=\cases{w_{\scr00(\a-1)}\dR(\xx_{(s-t-1)\b+1}\dR w_{\scr000})&for~$t=z$\cr
w_{\scr00(2\a-1)}\dR(w_{\scr00(\a-2)}\xx_{(s-t-1)\b+1}\dR\e)&for~$t=z+1$}=\e.}$$
{\it Subcase~\CaseIV.6:} $z+1<t<s$.
We find$$\eqalign{u_{ijk}&=(w_{\scr(t\a)0(\a-1)}
\dR w_{\scr(z\a)0(2\a-1)})\dR(w_{\scr(t\a)1(\a-1)}\xx_{(s-t-1)\b+1}
\dR w_{\scr(z\a)0(2\a-1)})\cr
&=w_{\scr((z-t-1)\a)0(2\a-1)}\dR\e=\e.}$$
{\it Subcase~\CaseIV.7:} $s<t<z$.
We find$$\eqalign{u_{ijk}&=(w_{\scr(s\a)1(\a-1)}\xx_{(t-s-1)\b+1}
\dR w_{\scr(s\a)1(\a-1)}\xx_{(z-s-1)\b+g})\dR(w_{\scr(s\a)0(\a-1)}
\dR w_{\scr(t\a)1(\a-1)}\xx_{(z-t-1)\b+g})\cr
&=(\xx_{(t-s-1)\b+1}\dR\xx_{(z-s-1)\b+g})\dR w_{\scr((t-s-1)\a)1(\a-1)}\xx_{(z-t-1)\b+g}=\e.}$$
{\it Subcase~\CaseIV.8:} $s<z\leq t\leq z+1$.
We find$$\eqalign{u_{ijk}&=(w_{\scr(s\a)1(\a-1)}\xx_{(t-s-1)\b+1}
\dR w_{\scr(s\a)1(\a-1)}\xx_{(z-s-1)\b+g}
)\dR(w_{\scr(s\a)0(\a-1)}\dR w_{\scr(t\a)1(z\a-t\a+\a)})\cr
&=(\xx_{(t-s-1)\b+1}
\dR\xx_{(z-s-1)\b+g}
)\dR w_{\scr((t-s-1)\a)1(z\a-t\a+\a)}=\e.}$$
{\it Subcase~\CaseIV.9:} $s+1<z+1<t$.
We find$$\eqalign{u_{ijk}&=(w_{\scr(s\a)1(\a-1)}\xx_{(t-s-1)\b+1}
\dR w_{\scr(s\a)1(\a-1)}\xx_{(z-s-1)\b+g})\dR(w_{\scr(s\a)0(\a-1)}
\dR w_{\scr(z\a)0(2\a-1)})\cr
&=(\xx_{(t-s-1)\b+1}
\dR\xx_{(z-s-1)\b+g})\dR w_{\scr((z-s-1)\a)0(2\a-1)}=\e.}$$
{\it Subcase~\CaseIV.10:} $z=s=t-1$.
We find$$\eqalign{u_{ijk}&=(w_{\scr(s\a)1(\a-1)}\xx_{(t-s-1)\b+1}
\dR w_{\scr(s\a)1(z\a-s\a+\a)}
)\dR(w_{\scr(s\a)0(\a-1)}
\dR w_{\scr(t\a)1(z\a-t\a+\a)})\cr
&=(w_{\scr(s\a)1(\a-1)}\xx_1\dR w_{\scr(s\a)1\a}
)\dR w_{\scr010}=(\xx_1\dR w_{\scr000}
)\dR w_{\scr010}=\e.}$$
{\it Subcase~\CaseIV.11:} $z=s<t-1$.
We find$$\eqalign{u_{ijk}&=(w_{\scr(s\a)1(\a-1)}\xx_{(t-s-1)\b+1}
\dR w_{\scr(s\a)1(z\a-s\a+\a)}
)\dR(w_{\scr(s\a)0(\a-1)}
\dR w_{\scr(z\a)0(2\a-1)})\cr
&=(w_{\scr(s\a)1(\a-1)}\xx_{(t-s-1)\b+1}
\dR w_{\scr(s\a)1\a}
)\dR w_{\scr00(\a-1)}=(\xx_{(t-s-1)\b+1}
\dR w_{\scr000})\dR w_{\scr00(\a-1)}=\e.}$$
{\it Subcase~\CaseIV.12:} $z=s+1<t-1$.
We find$$\eqalign{u_{ijk}&=(w_{\scr(s\a)1(\a-1)}\xx_{(t-s-1)\b+1}
\dR w_{\scr(s\a)1(\a-1)}\xx_{(z-s-1)\b+g}
)\dR(w_{\scr(s\a)0(\a-1)}
\dR w_{\scr(z\a)0(2\a-1)})\cr
&=(\xx_{(t-s-1)\b+1}\dR\xx_{g})\dR w_{\scr00(2\a-1)})=\e.}$$
{\it Subcase~\CaseIV.13:} $z+1<s<t$.
We find$$\eqalign{u_{ijk}&=(w_{\scr(s\a)1(\a-1)}\xx_{(t-s-1)\b+1}
\dR w_{\scr(z\a)0(2\a-1)}
)\dR(w_{\scr(s\a)0(\a-1)}
\dR w_{\scr(z\a)0(2\a-1)})\cr
&=w_{\scr((s-z-2)\a)1(\a-1)}\xx_{(t-s-1)\b+1}\dR\e=\e.}$$

\mgni{\it Case~\CaseV:} $e=1$,
$f\not=1$, $g\not=1$.

\ni{\it Subcase~\CaseV.1:} $s<z<t-1$.
We find$$\eqalign{u_{ijk}&=(w_{\scr(s\a)1(\a-1)}\xx_{(t-s-1)\b+f}
\dR w_{\scr(s\a)1(\a-1)}\xx_{(z-s-1)\b+g})\dR(w_{\scr(s\a)0(\a-1)}
\dR w_{\scr(z\a)0(2\a-1)})\cr
&=(\xx_{(t-s-1)\b+f}
\dR\xx_{(z-s-1)\b+g})\dR w_{\scr((z-s-1)\a)0(2\a-1)}=\e.}$$
{\it Subcase~\CaseV.2:} $s<t\leq z+1\leq t+1$.
We find$$\eqalign{u_{ijk}&=(w_{\scr(s\a)1(\a-1)}\xx_{(t-s-1)\b+f}
\dR w_{\scr(s\a)1(\a-1)}\xx_{(z-s-1)\b+g})\dR(w_{\scr(s\a)0(\a-1)}
\dR w_{\scr(t\a)f\a})\cr
&=(\xx_{(t-s-1)\b+f}\dR\xx_{(z-s-1)\b+g})\dR w_{\scr((t-s-1)\a)f\a}=\e.}$$
{\it Subcase~\CaseV.3:} $s<t=z-1$.
We find$$\eqalign{u_{ijk}&=(w_{\scr(s\a)1(\a-1)}\xx_{(t-s-1)\b+f}
\dR w_{\scr(s\a)1(\a-1)}\xx_{(z-s-1)\b+g})\cr
&\phantom{=5}\dR(w_{\scr(s\a)0(\a-1)}
\dR w_{\scr(t\a)f\a}w_{\scr00(\a-1)})\cr
&=(\xx_{(t-s-1)\b+f}\dR\xx_{(z-s-1)\b+g})\dR w_{\scr((t-s-1)\a)f\a}w_{\scr00(\a-1)}=\e.}$$
{\it Subcase~\CaseV.4:} $s<t<z-1$.
We find$$\eqalign{u_{ijk}&=(w_{\scr(s\a)1(\a-1)}\xx_{(t-s-1)\b+f}
\dR w_{\scr(s\a)1(\a-1)}\xx_{(z-s-1)\b+g})\cr
&\phantom{=5}\dR(w_{\scr(s\a)0(\a-1)}
\dR w_{\scr(t\a)f\a}w_{\scr00(\a-2)}\xx_{(z-t-2)\b+g}\cr
&=(\xx_{(t-s-1)\b+f}\dR\xx_{(z-s-1)\b+g})\dR w_{\scr((t-s-1)\a)f\a}w_{\scr00(\a-2)}\xx_{(z-t-2)\b+g}=\e.}$$
{\it Subcase~\CaseV.5:} $z\leq s\leq z+1<t$. We
find$$\eqalign{u_{ijk}&=(w_{\scr(s\a)1(\a-1)}\xx_{(t-s-1)\b+f}
\dR w_{\scr(s\a)1(z\a-s\a+\a)})\dR(w_{\scr(s\a)0(\a-1)}
\dR w_{\scr(z\a)0(2\a-1)})\cr
&=(w_{\scr(s\a)1(\a-1)}\xx_{(t-s-1)\b+f}
\dR w_{\scr(s\a)1(z\a-s\a+\a)})\dR\e=\e.}$$
{\it Subcase~\CaseV.6:} $t-1=s=z$.
We find$$\eqalign{u_{ijk}&=(w_{\scr(s\a)1(\a-1)}\xx_{(t-s-1)\b+f}
\dR w_{\scr(s\a)1(z\a-s\a+\a)})\dR(w_{\scr(s\a)0(\a-1)}
\dR w_{\scr(t\a)f\a})\cr
&=(\xx_f\dR w_{\scr000})\dR w_{\scr0f\a}=\e.}$$
{\it Subcase~\CaseV.7:} $z+1<s<t$.
We find$$\eqalign{u_{ijk}&=(w_{\scr(s\a)1(\a-1)}\xx_{(t-s-1)\b+f}
\dR w_{\scr(z\a)0(2\a-1)})\dR(w_{\scr(s\a)0(\a-1)}
\dR w_{\scr(z\a)0(2\a-1)})\cr
&=(w_{\scr(s\a)1(\a-1)}\xx_{(t-s-1)\b+f}
\dR w_{\scr(z\a)0(2\a-1)})\dR\e=\e.}$$
{\it Subcase~\CaseV.8:} $t\leq s\leq t+1<z$. We
find$$\eqalign{u_{ijk}&=(w_{\scr(s\a)1(t\a-s\a+\a)}
\dR w_{\scr(s\a)1(\a-1)}\xx_{(z-s-1)\b+g})\cr
&\phantom{=5}\dR(w_{\scr(t\a)f\a}
\dR w_{\scr(t\a)f\a}w_{\scr00(\a-2)}\xx_{(z-t-2)\b+g})\cr
&=\cases{w_{\scr00(\a-2)}\xx_{(z-s-1)\b+g}\dR w_{\scr00(\a-2)}\xx_{(z-t-2)\b+g}&for~$s=t+1$\cr
(w_{\scr000}
\dR\xx_{(z-s-1)\b+g})\dR w_{\scr00(\a-2)}\xx_{(z-t-2)\b+g}&for~$s=t$}=\e.}$$
{\it Subcase~\CaseV.9:} $s=t=z-1$.
We find$$\eqalign{u_{ijk}&=(w_{\scr(s\a)1(t\a-s\a+\a)}
\dR w_{\scr(s\a)1(\a-1)}\xx_{(z-s-1)\b+g})\dR(w_{\scr(t\a)f\a}
\dR w_{\scr(t\a)f\a}w_{\scr00(\a-1)})\cr
&=(w_{\scr000}
\dR\xx_g)\dR w_{\scr00(\a-1)}=\e.}$$
{\it Subcase~\CaseV.10.} Assume either~$s-1\leq t=z\leq s$ or~$s=t=z+1$.
We find$$u_{ijk}=(w_{\scr(s\a)1(t\a-s\a+\a)}
\dR w_{\scr(s\a)1(z\a-s\a+\a)})\dR(w_{\scr(t\a)f\a}
\dR w_{\scr(t\a)f\a})=\e.$$
{\it Subcase~\CaseV.11:} $t+1=z=s$.
We find$$\eqalign{u_{ijk}&=(w_{\scr(s\a)1(t\a-s\a+\a)}
\dR w_{\scr(s\a)1(z\a-s\a+\a)})\dR(w_{\scr(t\a)f\a}
\dR w_{\scr(t\a)f\a}w_{\scr00(\a-1)})=w_{\scr00(\a-1)}\dR w_{\scr00(\a-1)}=\e.}$$
{\it Subcase~\CaseV.12:} $z+1=t=s-1$.
We find$$u_{ijk}=(w_{\scr(s\a)1(t\a-s\a+\a)}
\dR w_{\scr(z\a)0(2\a-1)})\dR(w_{\scr(t\a)f\a}
\dR w_{\scr(t\a)f\a})=\e.$$
{\it Subcase~\CaseV.13:} $z+1<t\leq s\leq t+1$.
We find$$\eqalign{u_{ijk}&=(w_{\scr(s\a)1(t\a-s\a+\a)}
\dR w_{\scr(z\a)0(2\a-1)})\dR(w_{\scr(t\a)f\a}
\dR w_{\scr(z\a)0(2\a-1)})\cr
&=(w_{\scr(s\a)1(t\a-s\a+\a)}
\dR w_{\scr(z\a)0(2\a-1)})\dR\e=\e.}$$
{\it Subcase~\CaseV.14:} $t+1<s<z$.
We find$$\eqalign{u_{ijk}&=(w_{\scr(t\a)0(2\a-1)}
\dR w_{\scr(s\a)1(\a-1)}\xx_{(z-s-1)\b+g})\cr
&\phantom{=5}\dR(w_{\scr(t\a)f\a}w_{\scr00(\a-2)}\xx_{(s-t-2)\b+1}
\dR w_{\scr(t\a)f\a}w_{\scr00(\a-2)}\xx_{(z-t-2)\b+g})\cr
&=w_{\scr((s-t-2)\a)1(\a-1)}\xx_{(z-s-1)\b+g}\dR(\xx_{(s-t-2)\b+1}
\dR\xx_{(z-t-2)\b+g})=\e.}$$
{\it Subcase~\CaseV.15:} $t+1=z=s-1$.
We find$$\eqalign{u_{ijk}&=(w_{\scr(t\a)0(2\a-1)}
\dR w_{\scr(s\a)1(z\a-s\a+\a)})\dR(w_{\scr(t\a)f\a}w_{\scr00(\a-2)}\xx_{(s-t-2)\b+1}
\dR w_{\scr(t\a)f\a}w_{\scr00(\a-1)})\cr
&=w_{\scr010}\dR(\xx_1\dR w_{\scr000})=\e.}$$
{\it Subcase~\CaseV.16:} $t+1<z\leq s\leq z+1$.
We find$$\eqalign{u_{ijk}&=(w_{\scr(t\a)0(2\a-1)}
\dR w_{\scr(s\a)1(z\a-s\a+\a)})\dR(w_{\scr(t\a)f\a}w_{\scr00(\a-2)}\xx_{(s-t-2)\b+1}
\dR w_{\scr(t\a)f\a}w_{\scr00(\a-2)}\xx_{(z-t-2)\b+g})\cr
&=w_{\scr((s-t-2)\a)1(z\a-s\a+\a)}\dR(\xx_{(s-t-2)\b+1}\dR\xx_{(z-t-2)\b+g})=\e.}$$
{\it Subcase~\CaseV.17:} $z+1<t<s-1$.
We find$$u_{ijk}=(w_{\scr(t\a)0(2\a-1)}
\dR w_{\scr(z\a)0(2\a-1)})\dR(w_{\scr(t\a)f\a}w_{\scr00(\a-2)}\xx_{(s-t-2)\b+1}
\dR w_{\scr(z\a)0(2\a-1)})=\e\dR\e=\e.$$
{\it Subcase~\CaseV.18:} $t-1\leq z\leq t<s-1$.
We find$$u_{ijk}=(w_{\scr(t\a)0(2\a-1)}
\dR w_{\scr(z\a)0(2\a-1)})\dR(w_{\scr(t\a)f\a}w_{\scr00(\a-2)}\xx_{(s-t-2)\b+1}
\dR w_{\scr(t\a)f\a})=\e\dR\e=\e.$$
{\it Subcase~\CaseV.19:} $t+1=z<s-1$.
We find$$\eqalign{u_{ijk}&=(w_{\scr(t\a)0(2\a-1)}
\dR w_{\scr(z\a)0(2\a-1)})\dR(w_{\scr(t\a)f\a}w_{\scr00(\a-2)}\xx_{(s-t-2)\b+1}
\dR w_{\scr(t\a)f\a}w_{\scr00(\a-1)})\cr
&=w_{\scr00(\a-1)}\dR(\xx_{(s-t-2)\b+1}\dR w_{\scr000})=\e.}$$
{\it Subcase~\CaseV.20:} $t+1<z<s-1$.
We find$$\eqalign{u_{ijk}&=(w_{\scr(t\a)0(2\a-1)}
\dR w_{\scr(z\a)0(2\a-1)})\dR(w_{\scr(t\a)f\a}w_{\scr00(\a-2)}\xx_{(s-t-2)\b+1}
\dR w_{\scr(t\a)f\a}w_{\scr00(\a-2)}\xx_{(z-t-2)\b+g})\cr
&=w_{\scr((z-t-2)\a)0(2\a-1)}\dR(\xx_{(s-t-2)\b+1}\dR\xx_{(z-t-2)\b+g})=\e.}$$

\mgni{\it Case~\CaseVI:} $e\not=1$, $f=1$,
$g\not=1$.

\ni{\it Subcase~\CaseVI.1:} $t<z<s-1$.
We find$$\eqalign{u_{ijk}&=(w_{\scr(t\a)0(\a-1)}
\dR w_{\scr(z\a)0(2\a-1)}
)\dR(w_{\scr(t\a)1(\a-1)}\xx_{(s-t-1)\b+e}
\dR w_{\scr(t\a)1(\a-1)}\xx_{(z-t-1)\b+g})\cr
&=w_{\scr((z-t-1)\a)0(2\a-1)}\dR(\xx_{(s-t-1)\b+e}\dR\xx_{(z-t-1)\b+g})=\e.}$$
{\it Subcase~\CaseVI.2:} $z\leq t\leq z+1<s$.
We find$$\eqalign{u_{ijk}&=(w_{\scr(t\a)0(\a-1)}
\dR w_{\scr(z\a)0(2\a-1)}
)\dR(w_{\scr(t\a)1(\a-1)}\xx_{(s-t-1)\b+e}
\dR w_{\scr(t\a)1(z\a-t\a+\a)})\cr
&=\cases{w_{\scr00(\a-1)}\dR(\xx_{(s-t-1)\b+e}\dR w_{\scr000})&for~$t=z$\cr
\e\dR(w_{\scr00(\a-2)}\xx_{(s-t-1)\b+e}\dR\e)&for~$t=z+1$}=\e.}$$
{\it Subcase~\CaseVI.3:} $z+1<t<s$.
We find$$\eqalign{u_{ijk}&=(w_{\scr(t\a)0(\a-1)}
\dR w_{\scr(z\a)0(2\a-1)}
)\dR(w_{\scr(t\a)1(\a-1)}\xx_{(s-t-1)\b+e}
\dR w_{\scr(z\a)0(2\a-1)})=\e\dR\e=\e.}$$
{\it Subcase~\CaseVI.4:} $t<z\leq s\leq z+1$.
We find$$\eqalign{u_{ijk}&=(w_{\scr(t\a)0(\a-1)}
\dR w_{\scr(s\a)e\a}
)\dR(w_{\scr(t\a)1(\a-1)}\xx_{(s-t-1)\b+e}
\dR w_{\scr(t\a)1(\a-1)}\xx_{(z-t-1)\b+g})\cr
&=w_{\scr((s-t-1)\a)e\a}\dR(\xx_{(s-t-1)\b+e}\dR\xx_{(z-t-1)\b+g})=\e.}$$
{\it Subcase~\CaseVI.5:} $t<s=z-1$.
We find$$\eqalign{u_{ijk}&=(w_{\scr(t\a)0(\a-1)}
\dR w_{\scr(s\a)e\a}w_{\scr00(\a-1)}
)\dR(w_{\scr(t\a)1(\a-1)}\xx_{(s-t-1)\b+e}
\dR w_{\scr(t\a)1(\a-1)}\xx_{(z-t-1)\b+g})\cr
&=w_{\scr((s-t-1)\a)e\a}w_{\scr00(\a-1)}\dR(\xx_{(s-t-1)\b+e}\dR\xx_{(z-t-1)\b+g})=\e.}$$
{\it Subcase~\CaseVI.6:} $t<s<z-1$.
We find$$\eqalign{u_{ijk}&=(w_{\scr(t\a)0(\a-1)}
\dR w_{\scr(s\a)e\a}w_{\scr00(\a-2)}\xx_{(z-s-2)\b+g})\cr
&\phantom{=5}\dR(w_{\scr(t\a)1(\a-1)}\xx_{(s-t-1)\b+e}
\dR w_{\scr(t\a)1(\a-1)}\xx_{(z-t-1)\b+g})\cr
&=w_{\scr((s-t-1)\a)e\a}w_{\scr00(\a-2)}\xx_{(z-s-2)\b+g}\dR(\xx_{(s-t-1)\b+e}\dR\xx_{(z-t-1)\b+g})=\e.}$$
{\it Subcase~\CaseVI.7:} $z+1<s\leq t\leq
s+1$.
We find$$\eqalign{u_{ijk}&=(w_{\scr(s\a)e\a}
\dR w_{\scr(z\a)0(2\a-1)}
)\dR(w_{\scr(t\a)1(s\a-t\a+\a)}
\dR w_{\scr(z\a)0(2\a-1)})=\e\dR\e=\e.}$$
{\it Subcase~\CaseVI.8.} Assume either~$z+1=s=t$
or~$z=s=t$ or~$z=s=t-1$.
We find$$\eqalign{u_{ijk}&=(w_{\scr(s\a)e\a}
\dR w_{\scr(s\a)e\a}
)\dR(w_{\scr(t\a)1(s\a-t\a+\a)}
\dR w_{\scr(t\a)1(z\a-t\a+\a)})=\e\dR\e=\e.}$$
{\it Subcase~\CaseVI.9:} $z+1=s=t-1$.
We find$$\eqalign{u_{ijk}&=(w_{\scr(s\a)e\a}
\dR w_{\scr(s\a)e\a}
)\dR(w_{\scr(t\a)1(s\a-t\a+\a)}
\dR w_{\scr(z\a)0(2\a-1)})=\e\dR\e=\e.}$$
{\it Subcase~\CaseVI.10:} $s=t=z-1$.
We find$$\eqalign{u_{ijk}&=(w_{\scr(s\a)e\a}
\dR w_{\scr(s\a)e\a}w_{\scr00(\a-1)}
)\dR(w_{\scr(t\a)1(s\a-t\a+\a)}
\dR w_{\scr(t\a)1(\a-1)}\xx_{(z-t-1)\b+g})\cr
&=w_{\scr00(\a-1)}\dR(w_{\scr000}\dR\xx_g)=\e.}$$
{\it Subcase~\CaseVI.11:} $t=s+1=z$.
We find$$\eqalign{u_{ijk}&=(w_{\scr(s\a)e\a}
\dR w_{\scr(s\a)e\a}w_{\scr00(\a-1)}
)\dR(w_{\scr(t\a)1(s\a-t\a+\a)}
\dR w_{\scr(t\a)1(z\a-t\a+\a)})=w_{\scr00(\a-1)}\dR w_{\scr00(\a-1)}=\e.}$$
{\it Subcase~\CaseVI.12:} $z+1<s<t-1$.
We find$$\eqalign{u_{ijk}&=(w_{\scr(s\a)e\a}w_{\scr00(\a-2)}\xx_{(t-s-2)\b+1}
\dR w_{\scr(z\a)0(2\a-1)}
)\dR(w_{\scr(s\a)0(2\a-1)}
\dR w_{\scr(z\a)0(2\a-1)})=\e\dR\e=\e.}$$
{\it Subcase~\CaseVI.13:} $s-1\leq z\leq s<t-1$.
We find$$\eqalign{u_{ijk}&=(w_{\scr(s\a)e\a}w_{\scr00(\a-2)}\xx_{(t-s-2)\b+1}
\dR w_{\scr(s\a)e\a}
)\dR(w_{\scr(s\a)0(2\a-1)}
\dR w_{\scr(z\a)0(2\a-1)})=\e\dR\e=\e.}$$
{\it Subcase~\CaseVI.14:} $s+1=z=t-1$.
We find$$\eqalign{u_{ijk}&=(w_{\scr(s\a)e\a}w_{\scr00(\a-2)}\xx_{(t-s-2)\b+1}
\dR w_{\scr(s\a)e\a}w_{\scr00(\a-1)}
)\dR(w_{\scr(s\a)0(2\a-1)}
\dR w_{\scr(t\a)1(z\a-t\a+\a)})\cr
&=(\xx_1\dR w_{\scr000})\dR w_{\scr010}=\e.}$$
{\it Subcase~\CaseVI.15:} $s+1=z<t-1$.
We find$$\eqalign{u_{ijk}&=(w_{\scr(s\a)e\a}w_{\scr00(\a-2)}\xx_{(t-s-2)\b+1}
\dR w_{\scr(s\a)e\a}w_{\scr00(\a-1)}
)\dR(w_{\scr(s\a)0(2\a-1)}
\dR w_{\scr(z\a)0(2\a-1)})\cr
&=(\xx_{(t-s-2)\b+1}
\dR w_{\scr000}
)\dR w_{\scr00(\a-1)}=\e.}$$
{\it Subcase~\CaseVI.16:} $s+1<t<z$.
We find$$\eqalign{u_{ijk}&=(w_{\scr(s\a)e\a}w_{\scr00(\a-2)}\xx_{(t-s-2)\b+1}
\dR w_{\scr(s\a)e\a}w_{\scr00(\a-2)}\xx_{(z-s-2)\b+g})\cr
&\phantom{=5}\dR(w_{\scr(s\a)0(2\a-1)}
\dR w_{\scr(t\a)1(\a-1)}\xx_{(z-t-1)\b+g})\cr
&=(\xx_{(t-s-2)\b+1}\dR\xx_{(z-s-2)\b+g})\dR w_{\scr((t-s-2)\a)1(\a-1)}\xx_{(z-t-1)\b+g}=\e.}$$
{\it Subcase~\CaseVI.17:} $s+1<z\leq t\leq
z+1$.
We find$$\eqalign{u_{ijk}&=(w_{\scr(s\a)e\a}w_{\scr00(\a-2)}\xx_{(t-s-2)\b+1}
\dR w_{\scr(s\a)e\a}w_{\scr00(\a-2)}\xx_{(z-s-2)\b+g}
)\dR(w_{\scr(s\a)0(2\a-1)}
\dR w_{\scr(t\a)1(z\a-t\a+\a)})\cr
&=(\xx_{(t-s-2)\b+1}\dR\xx_{(z-s-2)\b+g})\dR w_{\scr((t-s-2)\a)1(z\a-t\a+\a)}=\e.}$$
{\it Subcase~\CaseVI.18:} $s+2<z+1<t$.
We find$$\eqalign{u_{ijk}&=(w_{\scr(s\a)e\a}w_{\scr00(\a-2)}\xx_{(t-s-2)\b+1}
\dR w_{\scr(s\a)e\a}w_{\scr00(\a-2)}\xx_{(z-s-2)\b+g}
)\dR(w_{\scr(s\a)0(2\a-1)}
\dR w_{\scr(z\a)0(2\a-1)})\cr
&=(\xx_{(t-s-2)\b+1}\dR\xx_{(z-s-2)\b+g})\dR w_{\scr((z-s-2)\a)0(2\a-1)})=\e.}$$

\mgni{\it Case~\CaseVII:} $e\not=1$,
$f\not=1$, $g=1$.

\ni{\it Subcase~\CaseVII.1:} $z<t<s-1$.
We find$$\eqalign{u_{ijk}&=(w_{\scr(t\a)0(2\a-1)}
\dR w_{\scr(z\a)0(\a-1)}
)\dR(w_{\scr(t\a)f\a}w_{\scr00(\a-2)}\xx_{(s-t-2)\b+e}
\dR w_{\scr(z\a)0(\a-1)})=\e\dR\e=\e.}$$
{\it Subcase~\CaseVII.2:} $t<s-1<z\leq s+1$.
We find$$\eqalign{u_{ijk}&=(w_{\scr(t\a)0(2\a-1)}
\dR w_{\scr(s\a)e\a}
)\dR(w_{\scr(t\a)f\a}w_{\scr00(\a-2)}\xx_{(s-t-2)\b+e}
\dR
w_{\scr(t\a)f\a}w_{\scr00(\a-2)}\xx_{(z-t-2)\b+1})\cr
&=w_{\scr((s-t-2)\a)e\a}
\dR(\xx_{(s-t-2)\b+e}\dR\xx_{(z-t-2)\b+1})=\e.}$$
{\it Subcase~\CaseVII.3:} $t+1<s<z-1$.
We find$$\eqalign{u_{ijk}&=(w_{\scr(t\a)0(2\a-1)}
\dR w_{\scr(s\a)e\a}w_{\scr00(\a-2)}\xx_{(z-s-2)\b+1})\cr
&\phantom{=5}\dR(w_{\scr(t\a)f\a}w_{\scr00(\a-2)}\xx_{(s-t-2)\b+e}
\dR
w_{\scr(t\a)f\a}w_{\scr00(\a-2)}\xx_{(z-t-2)\b+1})\cr
&=w_{\scr((s-t-2)\a)e\a}w_{\scr00(\a-2)}\xx_{(z-s-2)\b+1}
\dR(\xx_{(s-t-2)\b+e}\dR\xx_{(z-t-2)\b+1})}$$
{\it Subcase~\CaseVII.4:} $z<t=s-1$.
We find$$\eqalign{u_{ijk}&=(w_{\scr(s\a)e\a}
\dR w_{\scr(z\a)0(\a-1)}
)\dR(w_{\scr(t\a)f\a}w_{\scr00(\a-1)}
\dR w_{\scr(z\a)0(\a-1)})=\e\dR\e=\e.}$$
{\it Subcase~\CaseVII.5:} $z=t=s-1$.
We find$$\eqalign{u_{ijk}&=(w_{\scr(s\a)e\a}
\dR w_{\scr(z\a)0(\a-1)}
)\dR(w_{\scr(t\a)f\a}w_{\scr00(\a-1)}
\dR w_{\scr(t\a)f\a})=\e\dR\e=\e.}$$
{\it Subcase~\CaseVII.6:} $t+1=s=z$.
We find$$\eqalign{u_{ijk}&=(w_{\scr(s\a)e\a}
\dR w_{\scr(s\a)e\a}
)\dR(w_{\scr(t\a)f\a}w_{\scr00(\a-1)}
\dR w_{\scr(t\a)f\a})=\e\dR\e=\e.}$$
{\it Subcase~\CaseVII.7:} $t+1=s=z-1$.
We find$$\eqalign{u_{ijk}&=(w_{\scr(s\a)e\a}
\dR w_{\scr(s\a)e\a}
)\dR(w_{\scr(t\a)f\a}w_{\scr00(\a-1)}
\dR
w_{\scr(t\a)f\a}w_{\scr00(\a-2)}\xx_{(z-t-2)\b+1})=w_{\scr000}\dR\xx_1=\e.}$$
{\it Subcase~\CaseVII.8:} $t+1=s<z-1$.
We find$$\eqalign{u_{ijk}&=(w_{\scr(s\a)e\a}
\dR
w_{\scr(s\a)e\a}w_{\scr00(\a-2)}\xx_{(z-s-2)\b+1}
)\dR(w_{\scr(t\a)f\a}w_{\scr00(\a-1)}
\dR
w_{\scr(t\a)f\a}w_{\scr00(\a-2)}\xx_{(z-t-2)\b+1})\cr
&=w_{\scr00(\a-2)}\xx_{(z-s-2)\b+1}
\dR(w_{\scr000}\dR\xx_{(z-t-2)\b+1})=\e.}$$
{\it Subcase~\CaseVII.9:} $z<t=s$.
We find$$\eqalign{u_{ijk}&=(w_{\scr(s\a)e\a}
\dR w_{\scr(z\a)0(\a-1)}
)\dR(w_{\scr(t\a)f\a}
\dR w_{\scr(z\a)0(\a-1)})=\e\dR\e=\e.}$$
{\it Subcase~\CaseVII.10:} $s=t\leq z\leq t+1$.
We find$$\eqalign{u_{ijk}&=(w_{\scr(s\a)e\a}
\dR w_{\scr(s\a)e\a}
)\dR(w_{\scr(t\a)f\a}
\dR w_{\scr(t\a)f\a})=\e\dR\e=\e.}$$
{\it Subcase~\CaseVII.11:} $s=t<z-1$.
We find$$\eqalign{u_{ijk}&=(w_{\scr(s\a)e\a}
\dR
w_{\scr(s\a)e\a}w_{\scr00(\a-2)}\xx_{(z-s-2)\b+1}
)\dR(w_{\scr(t\a)f\a}
\dR
w_{\scr(t\a)f\a}w_{\scr00(\a-2)}\xx_{(z-t-2)\b+1})=\e.}$$
{\it Subcase~\CaseVII.12:} $z<s=t-1$.
We find$$\eqalign{u_{ijk}&=(w_{\scr(s\a)e\a}w_{\scr00(\a-1)}
\dR w_{\scr(z\a)0(\a-1)}
)\dR(w_{\scr(t\a)f\a}
\dR w_{\scr(z\a)0(\a-1)})=\e\dR\e=\e.}$$
{\it Subcase~\CaseVII.13:} $s=z=t-1$.
We find$$\eqalign{u_{ijk}&=(w_{\scr(s\a)e\a}w_{\scr00(\a-1)}
\dR w_{\scr(s\a)e\a}
)\dR(w_{\scr(t\a)f\a}
\dR w_{\scr(z\a)0(\a-1)})=\e\dR\e=\e.}$$
{\it Subcase~\CaseVII.14:} $s+1=t=z$.
We find$$\eqalign{u_{ijk}&=(w_{\scr(s\a)e\a}w_{\scr00(\a-1)}
\dR w_{\scr(s\a)e\a}
)\dR(w_{\scr(t\a)f\a}
\dR w_{\scr(t\a)f\a})=\e\dR\e=\e.}$$
{\it Subcase~\CaseVII.15:} $s+1=t=z-1$.
We find$$\eqalign{u_{ijk}&=(w_{\scr(s\a)e\a}w_{\scr00(\a-1)}
\dR
w_{\scr(s\a)e\a}w_{\scr00(\a-2)}\xx_{(z-s-2)\b+1}
)\dR(w_{\scr(t\a)f\a}
\dR w_{\scr(t\a)f\a})=w_{\scr000}\dR\xx_1=\e.}$$
{\it Subcase~\CaseVII.16:} $s+1=t<z-1$.
We find$$\eqalign{u_{ijk}&=(w_{\scr(s\a)e\a}w_{\scr00(\a-1)}
\dR
w_{\scr(s\a)e\a}w_{\scr00(\a-2)}\xx_{(z-s-2)\b+1}
)\dR(w_{\scr(t\a)f\a}
\dR
w_{\scr(t\a)f\a}w_{\scr00(\a-2)}\xx_{(z-t-2)\b+1})\cr
&=(w_{\scr000}\dR\xx_{(z-s-2)\b+1})\dR
w_{\scr00(\a-2)}\xx_{(z-t-2)\b+1}=\e.}$$
{\it Subcase~\CaseVII.17:} $z<s<t-1$.
We find$$\eqalign{u_{ijk}&=(w_{\scr(s\a)e\a}w_{\scr00(\a-2)}\xx_{(t-s-2)\b+f}
\dR w_{\scr(z\a)0(\a-1)}
)\dR(w_{\scr(s\a)0(2\a-1)}
\dR w_{\scr(z\a)0(\a-1)})=\e\dR\e=\e.}$$
{\it Subcase~\CaseVII.18:} $s\leq z\leq s+1<t$.
We find$$\eqalign{u_{ijk}&=(w_{\scr(s\a)e\a}w_{\scr00(\a-2)}\xx_{(t-s-2)\b+f}
\dR w_{\scr(s\a)e\a}
)\dR(w_{\scr(s\a)0(2\a-1)}
\dR w_{\scr(z\a)0(\a-1)})=\e\dR\e=\e.}$$
{\it Subcase~\CaseVII.19:} $s+1<z<t$.
We find$$\eqalign{u_{ijk}&=(w_{\scr(s\a)e\a}w_{\scr00(\a-2)}\xx_{(t-s-2)\b+f}
\dR
w_{\scr(s\a)e\a}w_{\scr00(\a-2)}\xx_{(z-s-2)\b+1}
)\dR(w_{\scr(s\a)0(2\a-1)}
\dR w_{\scr(z\a)0(\a-1)})\cr
&=(\xx_{(t-s-2)\b+f}\dR\xx_{(z-s-2)\b+1})\dR
w_{\scr((z-s-2)\a)0(\a-1)}=\e.}$$
{\it Subcase~\CaseVII.20:} $s+1<t\leq z\leq t+1$.
We find$$\eqalign{u_{ijk}&=(w_{\scr(s\a)e\a}w_{\scr00(\a-2)}\xx_{(t-s-2)\b+f}
\dR
w_{\scr(s\a)e\a}w_{\scr00(\a-2)}\xx_{(z-s-2)\b+1}
)\dR(w_{\scr(s\a)0(2\a-1)}
\dR w_{\scr(t\a)f\a})\cr
&=(\xx_{(t-s-2)\b+f}\dR\xx_{(z-s-2)\b+1})\dR
w_{\scr((t-s-2)\a)f\a}=\e.}$$
{\it Subcase~\CaseVII.21:} $s+1<t<z-1$.
We find$$\eqalign{u_{ijk}&=(w_{\scr(s\a)e\a}w_{\scr00(\a-2)}\xx_{(t-s-2)\b+f}
\dR w_{\scr(s\a)e\a}w_{\scr00(\a-2)}\xx_{(z-s-2)\b+1})\cr
&\phantom{=5}\dR(w_{\scr(s\a)0(2\a-1)}
\dR
w_{\scr(t\a)f\a}w_{\scr00(\a-2)}\xx_{(z-t-2)\b+1})\cr
&=(\xx_{(t-s-2)\b+f}\dR\xx_{(z-s-2)\b+1}
)\dR
w_{\scr((t-s-2)\a)f\a}w_{\scr00(\a-2)}\xx_{(z-t-2)\b+1}=\e.}$$

\mgni{\it Case~\CaseVIII:} $e\not=1$,
$f\not=1$, $g\not=1$.

\ni{\it Subcase~\CaseVIII.1:} $z+1<t<s-1$.
We find$$\eqalign{u_{ijk}&=(w_{\scr(t\a)0(2\a-1)}
\dR w_{\scr(z\a)0(2\a-1)}
)\dR(w_{\scr(t\a)f\a}w_{\scr00(\a-2)}\xx_{(s-t-2)\b+e}
\dR w_{\scr(z\a)0(2\a-1)})=\e\dR\e=\e.}$$
\ni{\it Subcase~\CaseVIII.2:} $t-1\leq z\leq t<s-1$.
We find$$\eqalign{u_{ijk}&=(w_{\scr(t\a)0(2\a-1)}
\dR w_{\scr(z\a)0(2\a-1)}
)\dR(w_{\scr(t\a)f\a}w_{\scr00(\a-2)}\xx_{(s-t-2)\b+e}
\dR w_{\scr(t\a)f\a})}=\e\dR\e=\e.$$
\ni{\it Subcase~\CaseVIII.3:} $z-1=t<s-1$.
We find$$\eqalign{u_{ijk}&=(w_{\scr(t\a)0(2\a-1)}
\dR w_{\scr(z\a)0(2\a-1)}
)\dR(w_{\scr(t\a)f\a}w_{\scr00(\a-2)}\xx_{(s-t-2)\b+e}
\dR w_{\scr(t\a)f\a}w_{\scr00(\a-1)})\cr
&=w_{\scr00(\a-1)}\dR(\xx_{(s-t-2)\b+e}\dR w_{\scr000})=\e.}$$
\ni{\it Subcase~\CaseVIII.4:} $t+1<z<s-1$.
We find$$\eqalign{u_{ijk}&=(w_{\scr(t\a)0(2\a-1)}
\dR w_{\scr(z\a)0(2\a-1)}
)\dR(w_{\scr(t\a)f\a}w_{\scr00(\a-2)}\xx_{(s-t-2)\b+e}
\dR w_{\scr(t\a)f\a}w_{\scr00(\a-2)}\xx_{(z-t-2)\b+g})\cr
&=w_{\scr((z-t-2)\a)0(2\a-1)}
\dR(\xx_{(s-t-2)\b+e}\dR\xx_{(z-t-2)\b+g})=\e.}$$
\ni{\it Subcase~\CaseVIII.5:} $t+1<z\leq s\leq z+1$.
We find$$\eqalign{u_{ijk}&=(w_{\scr(t\a)0(2\a-1)}
\dR w_{\scr(s\a)e\a}
)\dR(w_{\scr(t\a)f\a}w_{\scr00(\a-2)}\xx_{(s-t-2)\b+e}
\dR w_{\scr(t\a)f\a}w_{\scr00(\a-2)}\xx_{(z-t-2)\b+g})\cr
&=w_{\scr((s-t-2)\a)e\a}
\dR(\xx_{(s-t-2)\b+e}\dR\xx_{(z-t-2)\b+g})=\e.}$$
\ni{\it Subcase~\CaseVIII.6:} $t+1=z=s-1$.
We find$$\eqalign{u_{ijk}&=(w_{\scr(t\a)0(2\a-1)}
\dR w_{\scr(s\a)e\a}
)\dR(w_{\scr(t\a)f\a}w_{\scr00(\a-2)}\xx_{(s-t-2)\b+e}
\dR w_{\scr(t\a)f\a}w_{\scr00(\a-1)})\cr
&=w_{\scr0e\a}\dR(\xx_e\dR w_{\scr000})=\e.}$$
\ni{\it Subcase~\CaseVIII.7:} $t+1<s=z-1$.
We find$$\eqalign{u_{ijk}&=(w_{\scr(t\a)0(2\a-1)}
\dR w_{\scr(s\a)e\a}w_{\scr00(\a-1)}
)\cr
&\phantom{=5}\dR(w_{\scr(t\a)f\a}w_{\scr00(\a-2)}\xx_{(s-t-2)\b+e}
\dR w_{\scr(t\a)f\a}w_{\scr00(\a-2)}\xx_{(z-t-2)\b+g})\cr
&=w_{\scr((s-t-2)\a)e\a}w_{\scr00(\a-1)}
\dR(\xx_{(s-t-2)\b+e}\dR\xx_{(z-t-2)\b+g})=\e.}$$
\ni{\it Subcase~\CaseVIII.8:} $t+1<s<z-1$.
We find$$\eqalign{u_{ijk}&=(w_{\scr(t\a)0(2\a-1)}
\dR
w_{\scr(s\a)e\a}w_{\scr00(\a-2)}\xx_{(z-s-2)\b+g})\cr
&\phantom{=5}\dR(w_{\scr(t\a)f\a}w_{\scr00(\a-2)}\xx_{(s-t-2)\b+e}
\dR w_{\scr(t\a)f\a}w_{\scr00(\a-2)}\xx_{(z-t-2)\b+g})\cr
&=w_{\scr((s-t-2)\a)e\a}w_{\scr00(\a-2)}\xx_{(z-s-2)\b+g}
\dR(\xx_{(s-t-2)\b+e}\dR\xx_{(z-t-2)\b+g})=\e.}$$
\ni{\it Subcase~\CaseVIII.9:} $z+1<t=s-1$.
We find$$\eqalign{u_{ijk}&=(w_{\scr(s\a)e\a}
\dR w_{\scr(z\a)0(2\a-1)}
)\dR(w_{\scr(t\a)f\a}w_{\scr00(\a-1)}
\dR w_{\scr(z\a)0(2\a-1)})=\e\dR\e=\e.}$$
\ni{\it Subcase~\CaseVIII.10:} $z+1=t=s-1$.
We find$$\eqalign{u_{ijk}&=(w_{\scr(s\a)e\a}
\dR w_{\scr(z\a)0(2\a-1)}
)\dR(w_{\scr(t\a)f\a}w_{\scr00(\a-1)}
\dR w_{\scr(t\a)f\a})=\e\dR\e=\e.}$$
\ni{\it Subcase~\CaseVIII.11:} $z=t=s-1$.
We find$$\eqalign{u_{ijk}&=(w_{\scr(s\a)e\a}
\dR w_{\scr(s\a)e\a}
)\dR(w_{\scr(t\a)f\a}w_{\scr00(\a-1)}
\dR w_{\scr(t\a)f\a})=\e\dR\e=\e.}$$
\ni{\it Subcase~\CaseVIII.12:} $z-1=t=s-1$.
We find$$\eqalign{u_{ijk}&=(w_{\scr(s\a)e\a}
\dR w_{\scr(s\a)e\a}
)\dR(w_{\scr(t\a)f\a}w_{\scr00(\a-1)}
\dR w_{\scr(t\a)f\a}w_{\scr00(\a-1)})=\e\dR\e=\e.}$$
\ni{\it Subcase~\CaseVIII.13:} $t+1=s=z-1$.
We find$$\eqalign{u_{ijk}&=(w_{\scr(s\a)e\a}
\dR w_{\scr(s\a)e\a}w_{\scr00(\a-1)}
)\dR(w_{\scr(t\a)f\a}w_{\scr00(\a-1)}
\dR w_{\scr(t\a)f\a}w_{\scr00(\a-2)}\xx_{(z-t-2)\b+g})\cr
&=w_{\scr00(\a-1)}\dR(w_{\scr000}\dR\xx_g)=\e.}$$
\ni{\it Subcase~\CaseVIII.14:} $t+1=s<z-1$.
We find$$\eqalign{u_{ijk}&=(w_{\scr(s\a)e\a}
\dR w_{\scr(s\a)e\a}w_{\scr00(\a-2)}\xx_{(z-s-2)\b+g}
)\dR(w_{\scr(t\a)f\a}w_{\scr00(\a-1)}
\dR w_{\scr(t\a)f\a}w_{\scr00(\a-2)}\xx_{(z-t-2)\b+g})\cr
&=w_{\scr00(\a-2)}\xx_{(z-s-2)\b+g}
\dR(w_{\scr000}\dR\xx_{(z-t-2)\b+g})=\e.}$$
\ni{\it Subcase~\CaseVIII.15:} $z+1<s=t$.
We find$$\eqalign{u_{ijk}&=(w_{\scr(s\a)e\a}
\dR w_{\scr(z\a)0(2\a-1)}
)\dR(w_{\scr(t\a)f\a}
\dR w_{\scr(z\a)0(2\a-1)})=\e\dR\e=\e.}$$
\ni{\it Subcase~\CaseVIII.16:} $z\leq s=t\leq z+1$.
We find$$\eqalign{u_{ijk}&=(w_{\scr(s\a)e\a}
\dR w_{\scr(s\a)e\a}
)\dR(w_{\scr(t\a)f\a}
\dR w_{\scr(t\a)f\a})}=\e\dR\e=\e.$$
\ni{\it Subcase~\CaseVIII.17:} $s+1=t+1=z$.
We find$$\eqalign{u_{ijk}&=(w_{\scr(s\a)e\a}
\dR w_{\scr(s\a)e\a}w_{\scr00(\a-1)}
)\dR(w_{\scr(t\a)f\a}
\dR w_{\scr(t\a)f\a}w_{\scr00(\a-1)})=w_{\scr00(\a-1)}\dR
w_{\scr00(\a-1)}=\e.}$$
\ni{\it Subcase~\CaseVIII.18:} $s+1=t+1<z$.
We find$$\eqalign{u_{ijk}&=(w_{\scr(s\a)e\a}
\dR
w_{\scr(s\a)e\a}w_{\scr00(\a-2)}\xx_{(z-s-2)\b+g}
)\dR(w_{\scr(t\a)f\a}
\dR w_{\scr(t\a)f\a}w_{\scr00(\a-2)}\xx_{(z-t-2)\b+g})=\e.}$$
\ni{\it Subcase~\CaseVIII.19:} $z+1<s=t-1$.
We find$$\eqalign{u_{ijk}&=(w_{\scr(s\a)e\a}w_{\scr00(\a-1)}
\dR w_{\scr(z\a)0(2\a-1)})\dR(w_{\scr(t\a)f\a}
\dR w_{\scr(z\a)0(2\a-1)})=\e\dR\e=\e.}$$
\ni{\it Subcase~\CaseVIII.20:} $z+1=s=t-1$.
We find$$\eqalign{u_{ijk}&=(w_{\scr(s\a)e\a}w_{\scr00(\a-1)}
\dR w_{\scr(s\a)e\a})\dR(w_{\scr(t\a)f\a}
\dR w_{\scr(z\a)0(2\a-1)})=\e\dR\e=\e.}$$
\ni{\it Subcase~\CaseVIII.21:} $z=s=t-1$.
We find$$\eqalign{u_{ijk}&=(w_{\scr(s\a)e\a}w_{\scr00(\a-1)}
\dR w_{\scr(s\a)e\a}
)\dR(w_{\scr(t\a)f\a}
\dR w_{\scr(t\a)f\a})=\e\dR\e=\e.}$$
\ni{\it Subcase~\CaseVIII.22:} $z=s+1=t$.
We find$$\eqalign{u_{ijk}&=(w_{\scr(s\a)e\a}w_{\scr00(\a-1)}
\dR w_{\scr(s\a)e\a}w_{\scr00(\a-1)}
)\dR(w_{\scr(t\a)f\a}
\dR w_{\scr(t\a)f\a})=\e\dR\e=\e.}$$
\ni{\it Subcase~\CaseVIII.23:} $s+1=t=z-1$.
We find$$\eqalign{u_{ijk}&=(w_{\scr(s\a)e\a}w_{\scr00(\a-1)}
\dR
w_{\scr(s\a)e\a}w_{\scr00(\a-2)}\xx_{(z-s-2)\b+g}
)\dR(w_{\scr(t\a)f\a}
\dR w_{\scr(t\a)f\a}w_{\scr00(\a-1)})\cr
&=(w_{\scr000}
\dR\xx_g)\dR w_{\scr00(\a-1)}=\e.}$$
\ni{\it Subcase~\CaseVIII.24:} $s+1=t<z-1$.
We find$$\eqalign{u_{ijk}&=(w_{\scr(s\a)e\a}w_{\scr00(\a-1)}
\dR
w_{\scr(s\a)e\a}w_{\scr00(\a-2)}\xx_{(z-s-2)\b+g}
)\dR(w_{\scr(t\a)f\a}
\dR w_{\scr(t\a)f\a}w_{\scr00(\a-2)}\xx_{(z-t-2)\b+g})\cr
&=(w_{\scr000}\dR\xx_{(z-s-2)\b+g})
\dR w_{\scr00(\a-2)}\xx_{(z-t-2)\b+g}=\e.}$$
\ni{\it Subcase~\CaseVIII.25:} $z+1<s<t-1$.
We find$$\eqalign{u_{ijk}&=(w_{\scr(s\a)e\a}w_{\scr00(\a-2)}\xx_{(t-s-2)\b+f}
\dR w_{\scr(z\a)0(2\a-1)}
)\dR(w_{\scr(s\a)0(2\a-1)}
\dR w_{\scr(z\a)0(2\a-1)})=\e\dR\e.}$$
\ni{\it Subcase~\CaseVIII.26:} $s\leq z+1\leq s+1<t$.
We find$$\eqalign{u_{ijk}&=(w_{\scr(s\a)e\a}w_{\scr00(\a-2)}\xx_{(t-s-2)\b+f}
\dR w_{\scr(s\a)e\a}
)\dR(w_{\scr(s\a)0(2\a-1)}
\dR w_{\scr(z\a)0(2\a-1)})=\e\dR\e.}$$
\ni{\it Subcase~\CaseVIII.27:} $s+1=z<t-1$.
We find$$\eqalign{u_{ijk}&=(w_{\scr(s\a)e\a}w_{\scr00(\a-2)}\xx_{(t-s-2)\b+f}
\dR w_{\scr(s\a)e\a}w_{\scr00(\a-1)}
)\dR(w_{\scr(s\a)0(2\a-1)}
\dR w_{\scr(z\a)0(2\a-1)})\cr
&=(\xx_{(t-s-2)\b+f}\dR w_{\scr000})\dR
w_{\scr00(\a-1)}=\e.}$$
\ni{\it Subcase~\CaseVIII.28:} $s+1=z=t-1$.
We find$$\eqalign{u_{ijk}&=(w_{\scr(s\a)e\a}w_{\scr00(\a-2)}\xx_{(t-s-2)\b+f}
\dR w_{\scr(s\a)e\a}w_{\scr00(\a-1)}
)\dR(w_{\scr(s\a)0(2\a-1)}
\dR w_{\scr(t\a)f\a})\cr
&=(\xx_f\dR w_{\scr00(\a-1)})
\dR w_{\scr0f\a}=\e.}$$
\ni{\it Subcase~\CaseVIII.29:} $s+1<z<t-1$.
We find$$\eqalign{u_{ijk}&=(w_{\scr(s\a)e\a}w_{\scr00(\a-2)}\xx_{(t-s-2)\b+f}
\dR w_{\scr(s\a)e\a}w_{\scr00(\a-2)}\xx_{(z-s-2)\b+g}
)\dR(w_{\scr(s\a)0(2\a-1)}
\dR w_{\scr(z\a)0(2\a-1)})\cr
&=(\xx_{(t-s-2)\b+f}\dR\xx_{(z-s-2)\b+g}
)\dR w_{\scr((z-s-2)\a)0(2\a-1)}=\e.}$$
\ni{\it Subcase~\CaseVIII.30:} $s+1<z\leq t\leq z+1$.
We find$$\eqalign{u_{ijk}&=(w_{\scr(s\a)e\a}w_{\scr00(\a-2)}\xx_{(t-s-2)\b+f}
\dR w_{\scr(s\a)e\a}w_{\scr00(\a-2)}\xx_{(z-s-2)\b+g}
)\dR(w_{\scr(s\a)0(2\a-1)}
\dR w_{\scr(t\a)f\a})\cr
&=(\xx_{(t-s-2)\b+f}\dR\xx_{(z-s-2)\b+g})
\dR w_{\scr((t-s-2)\a)f\a}=\e.}$$
\ni{\it Subcase~\CaseVIII.31:} $s+1<t=z-1$.
We find$$\eqalign{u_{ijk}&=(w_{\scr(s\a)e\a}w_{\scr00(\a-2)}\xx_{(t-s-2)\b+f}
\dR w_{\scr(s\a)e\a}w_{\scr00(\a-2)}\xx_{(z-s-2)\b+g})\cr
&\phantom{=5}\dR(w_{\scr(s\a)0(2\a-1)}
\dR w_{\scr(t\a)f\a}w_{\scr00(\a-1)})\cr
&=\xx_{(t-s-2)\b+f}\dR\xx_{(z-s-2)\b+g})
\dR w_{\scr((t-s-2)\a)f\a}w_{\scr00(\a-1)}=\e.}$$
\ni{\it Subcase~\CaseVIII.32:} $s+1<t<z-1$.
We find$$\eqalign{u_{ijk}&=(w_{\scr(s\a)e\a}w_{\scr00(\a-2)}\xx_{(t-s-2)\b+f}
\dR w_{\scr(s\a)e\a}w_{\scr00(\a-2)}\xx_{(z-s-2)\b+g})\cr
&\phantom{=5}\dR(w_{\scr(s\a)0(2\a-1)}
\dR w_{\scr(t\a)f\a}w_{\scr00(\a-2)}\xx_{(z-t-2)\b+g})\cr
&=(\xx_{(t-s-2)\b+f}\dR\xx_{(z-s-2)\b+g}
)\dR
w_{\scr((t-s-2)\a)f\a}w_{\scr00(\a-2)}\xx_{(z-t-2)\b+g}=\e.}$$
This completes the proof of local coherence.

\vfill\eject
\bg\bgni\centerline{\scXII
References}

\bg\bg\bnni
\parindent=17pt

\Ref \\\Adj; S. I. Adjan; Defining relations and algorithmic
problems for groups and semi-groups; Proc. Steklov Inst. Math. {\bf
85} (1966).

\Reff \\\Art; E. Artin; Theorie der Z\"opfe; Abh.
Math. Sem. Hamburg Univ.; 4; 1925; 47--72.

\Ref \\\Bes;
D. Bessis; The dual braid monoid;
arXiv:math.GR/0101158.  

\Ref \\\BDM; D. Bessis, F. Digne
\& J. Michel; Springer theory in
braid groups and
the~Birman-Ko-Lee monoid;
arXiv:math.GR/0010254.

\Reff \\\Bir; J. Birman; Braids, links, and mapping class groups;
Annals of Math. Studies;  82; 1975;
Princeton Univ. Press, Princeton.

\Reff \\\BKL; J. Birman, K. H. Ko
\& S. J. Lee; A new approach to
the word and conjugacy
problems in the braid groups;
Advances in Math.;  139; 1998;
322--353.

\Ref \\\Bou; N. Bourbaki;
Groupes et alg\`ebres de Lie,
Chapitres 4, 5 et 6; Hermann,
Paris (1968).

\Reff \\\BMR; M. Brou\'e, G.
Malle \& R. Rouquier; Complex
reflection groups, braid groups,
Hecke algebras; J. Reine Angew.
Math.; 500; 1998; 127--190. 

\Reff \\\BrS; E. Brieskorn \& K.
Saito; Artin-Gruppen und
Coxeter-Gruppen; Invent. Math.;
17; 1972; 245--271.

\Ref \\\BuZ; G. Burde \& H. Zieschang; Knots; de
Gruyter, Berlin (1985).

\Ref \\\ClP; A. H. Clifford \& G.
B. Preston; The algebraic theory
of semigroups, vol.~1; AMS
Surveys {\bf 7} (1961).
 
\Reff \\\Dff; P. Dehornoy; Groups with a complemented
presentation; J. Pure Appl. Algebra; 116; 1997;
115--137.

\Ref \\\Dhh; ---; Braids and
Self-Distributivity; Progress in
Math. 192, Birkh\"auser (2000).

\Ref \\\Dgk; ---; Groupes
de Garside; Ann. Sci. Ec. Norm. Sup\'er. (2001) to appear.

\Ref \\\Djj; ---; Complete semigroup presentations; preprint
(2001).

\Reff \\\Dfx; P. Dehornoy \& L.
Paris; Gaussian groups and
Garside groups, two
generalizations of Artin groups;
Proc. London Math. Soc.; 79-3;
1999; 569--604.

\Reff \\\Del; P. Deligne; Les
immeubles des groupes de
tresses g\'en\'eralis\'es; Invent.
Math.; 17; 1972; 273--302.

\Ref \\\Eps; D. B. A. Epstein \& {\it al.};
Word Processing in Groups; Jones \& Bartlett
Publ. (1992).

\Reff \\\Gar; F. A. Garside; The
braid group and other groups;
Quart. J. Math. Oxford; 20; 1969;
235--254.

\Ref \\\HKo; J. W. Han \& K. H. Ko; Positive presentations of
the braid groups and the embedding problem;
arXiv:math.GR/0103132.

\Ref \\\Hig; P. J. Higgins; Techniques of semigroup theory;
Oxford Science Publications, Oxford University Press (1992)

\Ref \\\Hol; D. F. Holt; Automatic and hyperbolic groups;
to appear as Section B15 in "Handbook on the Heart of
Algebra" A. V. Mikhalev \& G. F. Pilz (eds.) Kluwer.

\Ref \\\LyS; R. C. Lyndon \& P. E.
Schupp; Combinatorial group theory; Springer-Verlag (1977).

\Ref \\\MKS; W. Magnus, A. Karrass \& D. Solitar;
Combinatorial group theory; Interscience, New York (1966).

\Reff \\\Pia; M. Picantin; The conjugacy problem in small
Gaussian groups; Comm. in~Algebra; 29-3; 2001; 1021--1039.

\Reff \\\Pib; ---; The
center of thin Gaussian groups; J. of Algebra; 145-1; 2001; 92-122.

\Ref \\\Pic; ---; Petits groupes
gaussiens; PhD Thesis,
Universit\'e de Caen (2000).

\Ref \\\Pid; ---; Diatomic Garside groups; in preparation.

\Reff \\\Rem; J. H. Remmers; On the geometry of semigroup
presentations; Advances in Math.; 36; 1980; 283--296.

\Ref \\\Rol; D. Rolfsen; Knots
and links; Publish or Perish, Inc
(1976).

\bye